\documentclass[12pt]{article}
\def\cdate{{November 22, 2026}}
\usepackage{amsfonts}
\usepackage{amsmath}
\usepackage{latexsym}
\usepackage{amssymb}
\usepackage{txfonts}
\usepackage[american]{babel}
\usepackage{graphicx}
\usepackage{framed}
\usepackage{xcolor}
\definecolor{mygray}{gray}{0.95} 
\definecolor{mydarkgray}{gray}{0.70} 
\colorlet{shadecolor}{mygray}

\makeatletter
\def\timenow{%
\@tempcnta=\time \divide\@tempcnta by 60 \number\@tempcnta:\multiply
\@tempcnta by 60 \@tempcntb=\time \advance\@tempcntb by -\@tempcnta
\ifnum\@tempcntb <10 0\number\@tempcntb\else\number\@tempcntb\fi}
\renewcommand{\@oddhead}{\small\it \@author\hfill\hfill\thepage}
\renewcommand{\@evenhead}{\small\it \@author\hfill\hfill\thepage}
\renewcommand{\@oddfoot}
{\vbox{\hrule\vspace{3pt}{\scriptsize\textit{\hfill\jobname.tex; \today; \timenow}}}}
\renewcommand{\@evenfoot}
{\vbox{\hrule\vspace{3pt}{\scriptsize\textit{\hfill\jobname.tex; \today; \timenow}}}}
\makeatother

\def\nmt{{
\null
\vspace{-4cm}
\par
\hspace*{50truemm}{\hrulefill}
\par
\vskip-4truemm
\par
\hspace*{50truemm}{\hrulefill}
\par\vskip5mm
\par
\hspace*{50truemm}{{\large\sc
New Mexico Tech {\rm
(\cdate)
}}}\vskip4mm
\par
\hspace*{50truemm}{\hrulefill}
\par
\vskip-4truemm
\par
\hspace*{50truemm}{\hrulefill}
\par}}

\def\RR{\mathbb{R}}
\def\CC{\mathbb{C}}

\def\cD{\mathcal{D}}
\def\cF{\mathcal{F}}

\def\cG{\mathcal{G}}
\def\cH{\mathcal{H}}

\def\cL{\mathcal{L}}
\def\cM{\mathcal{M}}

\def\cR{\mathcal{R}}

\def\cS{\mathcal{S}}
\def\cT{\mathcal{T}}

\def\cW{\mathcal{W}}

\def\slfrak{{\mathfrak{sl}}}

\def\ad{\mathrm{ad }}

\def\dim{\mathrm{dim\,}}

\def\Imag{\mathrm{Im\,}}

\def\Index{\mathrm{Index\,}}
\def\Ker{\mathrm{Ker\,}}
\def\Mat{\mathrm{Mat\,}}

\def\Ric{\mathrm{Ric}}

\def\Spin{\mathrm{Spin}}
\def\Span{\mathrm{Span}}

\def\tr{\mathrm{tr\,}}
\def\Tr{\mathrm{Tr\,}}
\def\vol{\mathrm{vol\,}}
\def\la{\langle}
\def\ra{\rangle}
\def\nn{{\nonumber}}

\def\proof{\par\noindent{\textbf{\textit{Proof:\;}}}}

\def\be{\begin{equation}}
\def\ee{\end{equation}}
\def\bea{\begin{eqnarray}}
\def\eea{\end{eqnarray}}
\def\bed{\begin{definition}{\ }}
\def\eed{\end{definition}}
\def\ed{\end{document}}
\def\bp{\begin{proposition}}
\def\ep{\end{proposition}}
\def\bc{\begin{center}}
\def\ec{\end{center}}
\def\bi{\begin{itemize}}
\def\ei{\end{itemize}}
\def\benum{\begin{enumerate}}
\def\eenum{\end{enumerate}}
\def\bmp{\begin{minipage}}
\def\emp{\end{minipage}}

\newtheorem{theorem}{Theorem}

\newtheorem{corollary}{Corollary}
\newtheorem{proposition}{Proposition}

\newtheorem{definition}{Definition}
\title{Geometry and Analysis of Spinor Valued Forms}
\author{Ivan G. Avramidi and Forest Good}

\begin{document}

\begin{titlepage}
\thispagestyle{empty}
\nmt

\bigskip
\bigskip
\bigskip
\bigskip
\centerline{\huge\bf Geometry and Analysis }
\bigskip
\bigskip
\centerline{\huge\bf of Spinor Valued Forms}
\bigskip
\bigskip
\bigskip
\centerline{\Large\bf Ivan G. Avramidi and Forest Good}
\bigskip
\centerline{\it New Mexico Institute of Mining and Technology}
\centerline{\it Socorro, NM 87801, USA}
\bigskip
\medskip

\begin{abstract}

We study spinor valued differential forms, that is, sections of the bundles 
$\cS_p=\cS\otimes\Lambda_p\otimes L$ of spinor valued $p$-forms twisted by a 
$U(1)$ line bundle $L$, on compact even-dimensional Riemannian spin 
(or $\Spin^\CC$) manifolds. 
We develop the algebra of spinor valued forms: the chirality, grading and 
duality involutions, the bilinear tensors built from spinor valued forms and 
the Fierz-type identities they satisfy, and the action of the Lie algebra 
$\slfrak(2,\RR)$ generated by Clifford multiplication and contraction, which 
we use to decompose the space of spinor valued forms into irreducible 
components and to construct an adapted orthonormal basis.
We introduce the twisted exterior derivative and coderivative, the Dirac and 
Rarita--Schwinger operators and, more generally, the natural first-order 
operators on spinor valued forms, as well as the twisted Hodge and signature 
operators, and derive their conformal transformation laws. 
We obtain the Weitzenb\"ock formulas for the twisted de Rham and Dirac 
Laplacians, express the curvature terms through the Weyl tensor, the Ricci 
tensor, the curvature of $L$ and the bilinear tensors of the spinor valued 
forms, and derive vanishing theorems, heat trace invariants and index formulas. 
Finally, in four dimensions we use these Weitzenb\"ock formulas to construct 
non-negative functionals whose absolute minima define two generalizations of 
the Seiberg--Witten monopole equations to right spinor valued $p$-forms, 
one based on the Dirac operator and one on the twisted de Rham Laplacian, 
and discuss the resulting integrability conditions.

\end{abstract}

\end{titlepage}


\section{Introduction}
\setcounter{equation}0

The Seiberg--Witten equations, motivated by the work of Seiberg and Witten 
on $\mathcal{N}=2$ supersymmetric Yang--Mills theory \cite{seiberg1994} and introduced 
in \cite{witten1994}, have had a
profound impact on four-manifold topology, providing invariants that are
both computable and geometrically meaningful. In the classical formulation,
one fixes a $\Spin^{\CC}$ structure on an oriented Riemannian $4$-manifold
$M$, considers the associated bundle $\cS^+ \otimes L$ of right spinors 
twisted by a $U(1)$ line bundle $L$, and studies a coupled system for a spinor
$\psi_+ \in C^\infty(\cS^+ \otimes L)$ and a unitary connection $A$ on $L$. 
The Seiberg--Witten invariants extracted from the moduli space of solutions have
been used to resolve long-standing questions in smooth $4$-manifold theory
\cite{witten1994, donaldson1996, morgan1996, moore2001, nicolaescu2000}.

The spinor fields appearing in the classical theory are sections of
$\cS^+ \otimes L$, which may be identified with spinor valued $0$-forms. It is
natural to ask whether an analogous theory can be developed for spinor valued
$p$-forms, that is, sections of
$\cS_p = \cS \otimes \Lambda_p \otimes L$, for all $p$. 
Such a generalization is not merely formal. 
Spinor valued $1$-forms are the Rarita--Schwinger fields 
\cite{rarita1941}, which describe the gravitino in supergravity 
\cite{freedman2012}. 
In Riemannian geometry, the Rarita--Schwinger operator and, more generally, 
the natural first-order operators (generalized gradients) on 
bundles associated with the spin group have been studied 
in \cite{branson1997, homma2019}. 
Spinor valued forms of all degrees are also closely related to 
the inhomogeneous differential forms of the Dirac--K\"ahler theory 
\cite{kahler1962}. 
Finally, gauge-theoretic equations involving spinors in higher
dimensions and on manifolds with special holonomy have been investigated
in \cite{tian2000, donaldson2010, haydys2012}.

The idea of replacing the spinor field in the Seiberg--Witten equations 
by a spinor valued form has been pursued recently for $p=1$. 
In \cite{nguyen22} a Seiberg--Witten type system built from the 
Rarita--Schwinger operator on spin $4$-manifolds was introduced and 
its $\mathrm{Pin}(2)$-equivariance was used to derive constraints on 
the topology of the manifold, under the assumption of compactness 
of the moduli space; the three-dimensional version of these equations 
was studied in \cite{haj23}, where a compactness theorem was proved. 
Our approach is different. We consider spinor valued forms of all degrees 
simultaneously, work with both the Dirac operator and the twisted 
exterior calculus on the whole space of spinor valued forms, and construct 
the monopole functionals directly from the Weitzenb\"ock formulas, 
which makes the role of the curvature, in particular, the Weyl tensor, 
explicit.

In this paper we develop the general geometric and analytic framework 
for such a generalization on even-dimensional Riemannian manifolds. 
Let $M$ be a compact, oriented, $n$-dimensional Riemannian manifold, 
$n=2\nu$, with a $\Spin^{\CC}$ structure 
(every oriented $4$-manifold has one, but this is not true in 
higher dimensions). 
If $M$ is spin, let $\cS$ be the spinor bundle and $L$ be a Hermitian 
line bundle over $M$. In general, neither $\cS$ nor $L$ is defined 
globally, but the twisted bundle $\cS\otimes L$ is globally defined 
(with $L$ playing the role of the square root of the determinant 
line bundle of the $\Spin^\CC$ structure). 
Since all our constructions involve only the bundle $\cS\otimes L$, 
they apply in the general case without modification. 
Let $A$ be a unitary connection $1$-form on $L$ 
with the curvature $2$-form $F = dA$.

For each $p = 0, 1, \ldots, n$, we define the twisted bundle
\begin{equation}
  \cS_p = \cS \otimes \Lambda_p \otimes L,
\end{equation}
where $\Lambda_p = \Lambda_p T^*M$. The classical Seiberg--Witten theory
corresponds to the case $p = 0$ in four dimensions.

Our main results are the following.

\benum
\item
{\it Algebra.} 
We study the bilinear tensors constructed from spinor valued forms 
and the Fierz-type identities they satisfy; in particular, 
we show that the even-degree bilinear forms of chiral spinors satisfy 
the duality relations $D^{\pm *}_{(2k)}=\pm D^\pm_{(2\nu-2k)}$. 
In four dimensions we decompose the spaces of right and left 
spinor valued forms into irreducible representations of 
the group $\Spin(4)=SU(2)\times SU(2)$. 
We show that the operators $\alpha$ and $\beta$ of exterior and 
interior Clifford multiplication, together with the degree operator, 
generate a representation of the Lie algebra $\slfrak(2,\RR)$ 
on the space of spinor valued forms, analogous to the Lefschetz 
$\slfrak(2)$ action in K\"ahler geometry. 
We determine the kernels of these operators and show that 
the subspace annihilated by both $\alpha$ and $\beta$ is 
trivial except in the middle degree, $p=\nu$. 
We decompose the space of spinor valued forms into irreducible 
$\slfrak(2,\RR)$ components and construct an adapted orthonormal basis.

\item
{\it Differential operators.} 
We introduce the twisted exterior derivative $D$ and coderivative $\bar D$ 
and the natural first-order operators $\cT_{p,k}$ on spinor valued 
$p$-forms; this family includes the Dirac operator $\cD_p$ and 
the Rarita--Schwinger operator $\cR_p$, which are the only non-zero 
ones in four dimensions. We show that the Dirac operator 
is expressed in terms of the operators $D$ and $\bar D$ 
and vice versa, and derive the conformal transformation laws 
of these operators. We introduce the operators 
$\cL_1=D+\bar D$ and $\cL_2=i(D-\bar D)$ on the total space 
of spinor valued forms and the duality 
involutions $\tau_1$ and $\tau_2$ anticommuting with them; 
the restrictions of $\cL_j$ to even and odd forms and to self-dual and 
anti-self-dual forms define the twisted Hodge operators and the twisted 
signature operators.

\item
{\it Weitzenb\"ock formulas.} 
We derive the Weitzenb\"ock formulas for the twisted de Rham Laplacian 
$L_p$ and for the square of the Dirac operator, 
\be
L_p=-\Delta+W_p
\ee
 and 
\be
-\cD^2_p=-\Delta+V_p,
\ee 
and express the quadratic forms $\la\psi,W_p\psi\ra$ and 
$\la\psi,V_p\psi\ra$ in terms of the Weyl tensor, 
the Ricci tensor, the scalar curvature, the curvature $F$ 
and the bilinear tensors of $\psi$. 
This gives vanishing theorems for twisted harmonic forms 
and for the zero modes of the Dirac operator.

\item
{\it Global invariants.}
We study the heat traces of these operators and express 
the indices of the Dirac operator, the twisted Hodge operators and 
the twisted signature operators in terms of the heat kernel 
coefficients \cite{gilkey1995, berline2004, avramidi2015}.

\item
{\it Monopoles.} 
In four dimensions, we use the Weitzenb\"ock formulas 
to construct two non-negative functionals on right spinor 
valued $p$-forms by completing the squares. 
The first one, $S_D$, is built from the Dirac operator and reduces 
to the Seiberg--Witten functional for $p=0$; 
its absolute minima are given by the equations 
\be
\cD\psi_+=0, \qquad
F_+=-iB_+(\psi).
\ee 
where the $2$-form $B_+(\psi)$, defined in (\ref{929igax}), is the natural 
generalization of the Seiberg--Witten moment map.
We show that, unlike the case $p=0$, for $p\ge 1$ 
the quartic term $|B_+|^2$ is not coercive, since it vanishes 
when the matrix $\Phi=\frac{1}{p!}\psi_+\otimes\bar\psi_+$ is 
proportional to the identity. 
The second one, $S_{dR}$, is built from the twisted de Rham Laplacian; 
its absolute minima are given by the system 
\be
D\psi_+=\bar D\psi_+=0, 
\qquad
F=-i\Omega(\psi),
\ee 
where $\Omega(\psi)$ is the antisymmetric part of the bilinear tensor
built from $\psi_+$ defined by (\ref{98igaxz}).
This system is overdetermined
and implies the algebraic integrability constraints 
$Q_p\psi_+=\bar Q_p\psi_+=0$ and $d\Omega=0$.
\eenum

The paper is organized as follows. 
In Section~\ref{sec:clifford} we review the basic properties
of the Clifford algebra and introduce the notation.
In Section~\ref{sec:forms} we introduce spinor valued forms,
the associated bilinear tensors and their Fierz identities, 
and discuss the representation theory in four dimensions.
In Section~\ref{sec:algebra} we introduce the $\slfrak(2,\RR)$ algebra 
of operators on spinor valued forms and study its kernels and images.
In Section~\ref{sec:basis} we construct an orthonormal basis
of spinor valued forms adapted to this algebra. 
In Section~\ref{sec:geometry} we review the relevant Riemannian
geometry: conformal transformations, the decomposition of 
the Weyl tensor in four dimensions, and the topological 
restrictions on Einstein metrics.
In Section~\ref{sec:calculus} we develop the twisted exterior
calculus, construct the Dirac, Rarita--Schwinger and 
twisted Hodge operators and the corresponding Laplacians, 
and derive the Weitzenb\"ock formulas.
In Section~\ref{sec:analysis} we analyze the curvature terms 
in the Weitzenb\"ock formulas and the resulting vanishing theorems, 
and discuss the heat traces and the indices of these operators. 
In Section~\ref{sec:monopoles} we review the Seiberg--Witten 
monopoles and propose the generalized monopole equations.


\section{Clifford Algebra}
\setcounter{equation}0
\label{sec:clifford}


Let $\Mat(2^\nu,\CC)$ be the set of complex $2^{\nu}\times 2^{\nu}$ matrices.
The Dirac matrices $\gamma_a\in \Mat(2^\nu,\CC)$, $a=1,\dots, 2\nu,$ 
are traceless involutive complex matrices
generating the Clifford algebra
\be
\gamma_a\gamma_b+\gamma_b\gamma_a=2\delta_{ab}I,
\label{12igax}
\ee
where $I$ is the $2^{\nu}$-dimensional unit matrix.
It is easy to see that
the Hermitian conjugate of the Dirac matrices also generate the Clifford algebra;
therefore, there must exist a matrix $\eta$ such that
\be
\gamma_a^*=\eta\gamma_a\eta^{-1},
\ee
which can always be chosen to be Hermitian, $\eta^*=\eta$.
One can show that there is a representation in which all Dirac matrices are Hermitian,
$\gamma_a^*=\gamma_a$,
and, therefore, one can choose the matrix $\eta$
to be equal to the unit matrix, $\eta=I$;
everywhere below we choose this matrix.
Then the Dirac matrices are traceless Hermitian involutions,
\be
(\gamma_a)^2=I, \qquad
\gamma_a^*=\gamma_a,\qquad
\tr\gamma_a=0.
\ee
Further, we define the anti-symmetrized
products of the Dirac matrices
\be
\gamma_{a_1\dots a_k}=\gamma_{[a_1}\cdots \gamma_{a_k]},
\qquad
k=1,\dots,2\nu.
\ee
Of course, these matrices vanish for $k>2\nu$.
For $ 1\le k\le 2\nu $,
all these matrices
are traceless and are (anti)-Hermitian (anti)-involutions,
\be
(\gamma_{a_1\dots a_k})^2 = (-1)^{k(k-1)/2}I,
\qquad
\tr \gamma_{a_1\dots a_k}=0.
\qquad
\gamma_{a_1\dots a_k}^{*} = (-1)^{k(k-1)/2}
{}\gamma_{a_1\dots a_k}{}.
\ee
Therefore,
\be
\gamma_{a_1\dots a_k}\gamma^{a_1\dots a_k} 
= (-1)^{k(k-1)/2}\frac{(2\nu)!}{(2\nu-k)!}I,
\ee

The traces of the products of these matrices are
\bea
\tr  \gamma_{a_1\cdots a_k}\gamma^{b_1\cdots b_m} &=& 0,
\qquad k\ne m,
\\[5pt]
\tr  \gamma_{a_1\cdots a_k}\gamma^{b_1\cdots b_k} &=&
(-1)^{k(k-1)/2}
\delta^{b_1\dots b_k}_{a_1\dots a_k}2^\nu,
\eea
where
\be
\delta^{b_1\dots b_k}_{a_1\dots a_k}
=k!\delta^{[b_1}_{[a_1}\cdots \delta^{b_k]}_{a_k]}
\ee
are the generalized Kronecker symbols.

The products of these matrices can be expressed as a linear combination
of them \cite{zhelnorovich2019},
\be
\gamma_{a_1\cdots a_k}\gamma^{b_1\cdots b_m}
=\sum_{p=0}^{\min{\{m,k\}}}
(-1)^{p(2k-p-1)/2}{k\choose p}{m\choose p}
\delta^{[b_1\dots b_p}_{[a_1\dots a_p}\gamma_{a_{p+1}\cdots a_k]}
{}^{b_{p+1}\cdots b_m]},
\ee
In this formula, in the case when $k=m$ there is a term with $p=k=m$;
the gamma matrix without indices is defined to be the
unit matrix. More precisely,
\bea
\gamma_{a_1\cdots a_k}\gamma^{b_1\cdots b_k}
&=&\sum_{p=0}^{k-1}
(-1)^{p(2k-p-1)/2}{k\choose p}^2
\delta^{[b_1\dots b_p}_{[a_1\dots a_p}\gamma_{a_{p+1}\cdots a_k]}
{}^{b_{p+1}\cdots b_k]}
\nn\\
&&+(-1)^{k(k-1)/2}
\delta^{b_1\dots b_k}_{a_1\dots a_k}I.
\eea

For $k =1, 2$ we have
\bea
\gamma_{a}\gamma^{b_1\cdots b_m}
&=&
\gamma_{a}{}^{b_{1}\cdots b_m}
+m\delta^{[b_1}_{a}\gamma^{b_{2}\cdots b_m]},
\\
\gamma_{a_1a_2}\gamma^{b_1\cdots b_m}
&=&\gamma_{a_1a_2}{}^{b_1\cdots b_m}
-2m\delta_{[a_1}^{[b_1}\gamma_{a_2]}{}^{b_2\cdots b_m]}
-m(m-1)\delta_{[a_1}^{[b_1}\delta_{a_2]}^{b_2}\gamma^{b_3\cdots b_m]},
\eea
in particular,
\bea
\gamma_{c}\gamma^{ab}
&=&\gamma_{c}{}^{ab}
+2\delta^{[a}_{c}\gamma^{b]},
\\
\gamma_{cd}\gamma^{ab}
&=&\gamma_{cd}{}^{ab}
-4\delta^{[a}_{[c}\gamma_{d]}{}^{b]}
-2\delta^{[a}_{[c}\delta^{b]}_{d]},
\label{gamma2}
\eea
We also have the commutators
\bea
[\gamma_{a},\gamma^{b_1\cdots b_m}]
&=&
[1-(-1)^m]\gamma_{a}{}^{b_{1}\cdots b_m}
+m[1+(-1)^m]\delta^{[b_1}_{a}\gamma^{b_{2}\cdots b_m]},
\\{}
[\gamma_{a_1a_2},\gamma^{b_1\cdots b_m}]
&=&
-4m\delta_{[a_1}^{[b_1}\gamma_{a_2]}{}^{b_2\cdots b_m]}
\eea
Therefore, the matrices $\gamma_{ab}$ form a representation
of the spin group
\be
[\gamma_{cd},\gamma^{ab}]
=-8\delta^{[a}_{[c}\gamma_{d]}{}^{b]}.
\ee

The collection of matrices
\be
I, \ \gamma_{a_1}, \ \gamma_{a_1a_2}, \ \cdots , \ \gamma_{a_1\dots a_{2\nu}}
\ee
with $a_1<\cdots<a_{2\nu}$, is linearly independent and forms a basis
in the space $\Mat(2^\nu,\CC)$. Therefore, any such matrix
$M\in\Mat(2^\nu,\CC)$ can be written in the form
\be
M=2^{-\nu}\left\{
C I+ \sum_{k=1}^{2\nu}\frac{1}{k!}(-1)^{k(k-1)/2}
C_{a_1\dots a_k}\gamma^{a_1\dots a_k}
\right\},
\ee
where
\bea
C &=& \tr M,
\\
C_{a_1\dots a_k} &=& \tr M\gamma_{a_1\dots a_k}.
\eea

Further, we define the chirality matrix
\be
\Gamma{}
=i^{\nu}\gamma_1\cdots\gamma_{2\nu},
\ee
so that
\be
\gamma_{a_1\dots a_{2\nu}}
=i^{-\nu}\varepsilon_{a_1\dots a_{2\nu}}
\Gamma,
\ee
This matrix is also traceless Hermitian involution 
\be
\Gamma^2 =I,
\qquad
\Gamma^*=\Gamma,
\qquad
\tr\Gamma=0,
\ee
which anti-commutes with the Dirac matrices,
\be
\Gamma{}\gamma_{a}=-\gamma_{a}\Gamma{},
\ee
and, more generally,
\be
\Gamma{}\gamma_{a_1\dots a_k}=(-1)^k\gamma_{a_1\dots a_k}\Gamma{}.
\ee

The matrices
\be
P_\pm=\frac{1}{2}(I\pm \Gamma{})
\ee
define complementary orthogonal projections.
and, therefore,
\bea
P_+ \gamma_{a_1\dots a_{2k+1}}P_+ &=& P_-\gamma_{a_1\dots a_{2k+1}}P_- = 0,
\label{213iga}
\\
P_+ \gamma_{a_1\dots a_{2k}}P_- &=&
P_- \gamma_{a_1\dots a_{2k}}P_+ = 0.
\eea
This motivates the definition of the chiral part of these matrices
\bea
\gamma^\pm_{a_1\dots a_{2k}}=P_\pm \gamma_{a_1\dots a_{2k}}P_\pm.
\eea

The dual matrices
\be
\tilde \gamma_{a_1\dots a_{k}}=\frac{1}{(2\nu-k)!}
\varepsilon_{b_1\dots b_{2\nu-k} a_1\dots a_{k}}\gamma^{b_1\dots b_{2\nu-k}},
\ee
are (anti)-Hermitian,
\bea
\tilde\gamma_{a_1\dots a_k}^* &=& (-1)^{\nu+k(k+1)/2}
{}\tilde\gamma_{a_1\dots a_k}{}.
\eea
and have the form
\bea
\tilde\gamma_{a_1\dots a_k}
&=& (-1)^{k(k-1)/2}i^{-\nu}
\Gamma\gamma_{a_1\dots a_k}.
\eea

Most of the time we will
restrict to four dimensions, $n=2\nu=4$.
Let $\sigma^i$, $i=1,2,3$, be the Pauli matrices
\bea
\sigma_1=
\left(
\begin{array}{cc}
0 & 1\\
1 & 0\\
\end{array}
\right),
\qquad
\sigma_2=
\left(
\begin{array}{cc}
0 & -i\\
i & 0\\
\end{array}
\right),
\qquad
\sigma_3=
\left(
\begin{array}{cc}
1 & 0\\
0 & -1\\
\end{array}
\right)\,.
\eea
These are Hermitian traceless involutive matrices
with the products
\be
\sigma_i\sigma_j=\delta_{ij}I+i\varepsilon_{ijk}\sigma_k\,.
\ee
They satisfy the following commutation relations
\be
[\sigma_i,\sigma_j]=2i\varepsilon_{ijk}\sigma_k
\ee
and generate two-dimensional Clifford algebra
\be
\sigma_i\sigma_j+\sigma_j\sigma_i=2\delta_{ij}I.
\ee
The dual matrices are defined by
\be
\tilde\sigma_{jk}=\varepsilon_{ljk}\sigma^l.
\ee

The Dirac matrices can be chosen in the form
(here $j=1,2,3$)
\be
\gamma^j
=
\left(
\begin{array}{cc}
0 & -i\sigma^j\\
i\sigma^j & 0\\
\end{array}
\right)\,,
\qquad
\gamma^4
=\left(\begin{array}{cc}
0 & I\\
I & 0\\
\end{array}\right)
\,,
\ee
The chirality matrix and the chiral projections are
\be
\Gamma{}{}=-\gamma_{1234}
=
\left(\begin{array}{cc}
	-I & 0\\
	0 & I\\
\end{array}\right);
\ee
\be
P_-
=\left(\begin{array}{cc}
	I & 0\\
	0 & 0\\
\end{array}\right), \qquad
P_+
=
\left(\begin{array}{cc}
	0 & 0\\
	0 & I\\
\end{array}\right).
\ee
Next, we have
\be
\gamma_{abc} = -\varepsilon_{abcd}\Gamma\gamma^d,
\ee
and, of course,
\be
\gamma_{abcd} = -\varepsilon_{abcd}\Gamma.
\ee

The matrices $\gamma_{ab}$ are block-diagonal
(here $j,k=1,2,3$)
\bea
\gamma_{jk} =
\left(\begin{array}{cc}
i\tilde\sigma_{jk} & 0\\
0 & i\tilde\sigma_{jk}\\
\end{array}\right),
\qquad
\gamma_{j4} =
\left(\begin{array}{cc}
-i\sigma_j & 0\\
0 & i\sigma_j\\
\end{array}\right).
\eea
The dual matrices
\be
\tilde\gamma_{ab} = \frac{1}{2}\varepsilon_{cdab}\gamma^{cd}
\ee
have the form
(here $j,k=1,2,3$)
\bea
\tilde\gamma_{jk} &=&
%
\left(\begin{array}{cc}
	-i\tilde\sigma_{jk} & 0\\
	0 & i\tilde\sigma_{jk}\\
\end{array}\right),
\qquad
\tilde\gamma_{j4} =
\left(\begin{array}{cc}
	i\sigma_j & 0\\
	0 & i\sigma_j\\
\end{array}\right).
\eea
Therefore,
\be
\tilde\gamma_{ab}
=\Gamma{}{}\gamma_{ab},
\ee

The chiral parts of these matrices 
\bea
\gamma^-_{jk} &=& P_-\gamma_{jk}=
\left(\begin{array}{cc}
	i\tilde\sigma_{jk} & 0\\
	0 & 0\\
\end{array}\right),
\qquad
\gamma^-_{j4} =P_-\gamma_{j4}=
\left(\begin{array}{cc}
	-i\sigma_j & 0\\
	0 & 0\\
\end{array}\right),
\\
\gamma^+_{jk} &=&
P_+\gamma_{jk}=
\left(\begin{array}{cc}
	0 & 0\\
	0 & i\tilde\sigma_{jk}\\
\end{array}\right),
\qquad
\gamma^+_{j4} = P_+\gamma_{j4}=
\left(\begin{array}{cc}
	0 & 0\\
	0 & i\sigma_j\\
\end{array}\right).
\eea
are 
anti-self-dual and self-dual, that is,
\be
\gamma^\pm_{ab}
=P_{\pm}\gamma_{ab}
=\frac{1}{2}\left(\gamma_{ab}
\pm\tilde\gamma_{ab}\right).
\ee

The product of the chiral parts of the matrices $\gamma_{ab}$ we get
\be
\gamma^\pm_{cd}\gamma_\pm^{ab}
=4\delta^{[a}{}_{[c}\gamma_\pm^{b]}{}_{d]}
-\left(2\delta^{[a}{}_{[c}\delta^{b]}{}_{d]}
\pm\varepsilon_{cd}{}^{ab}
\right)P_\pm,
\label{gamma2x}
\ee
and $\gamma^+_{ab}\gamma^-_{cd}=0$.
We will need the identities
\bea
\gamma_{cb}\gamma^{ab}
&=&-2\gamma_c{}^a-3\delta^a_c,
\\
\gamma_{ab}\gamma^{ab} &=& -12,
\eea
We compute some important products
\bea
\gamma_{cd}\gamma^{ab}\gamma_e{}^{d}
&=&
2 \gamma^{a b}\,_{ce}\,\discretionary{}{}{}
+ \delta_{ce}\,\gamma^{a b}\discretionary{}{}{}
+4 
\delta^{[a}_e \delta^{b]}_{c},
\\
\gamma_{cd}\gamma^{ab}\gamma_b{}^{d}&=&
\gamma^{a}\,_{c}
-6\delta^{a}_c,
\\
\gamma_{ad}\gamma^{ab}\gamma_b{}^{d}&=&
-24.
\eea
We will also need the identity
\bea
\gamma^b\gamma_{ba_1\dots a_{k}}
&=& \gamma_{a_1\dots a_{k}b}\gamma^{b}
=(4-k)\gamma_{a_{1}\cdots a_{k}},
\eea


\section{Spinor Valued Forms}
\setcounter{equation}0
\label{sec:forms}

Let $S$ be the complex vector space of spinors;
its complex dimension is
$\dim S=2^\nu$.
We define the conjugate spinor $\bar\psi$,
\be
\bar\psi=\psi^*{},
\ee
the inner product and the norm,
\be
\langle \chi,\psi\rangle=\bar\chi\psi,
\qquad |\psi|^2=\langle\psi,\psi\rangle.
\ee

Every spinor $\psi$ defines a collection of real valued forms by bilinear products,
\bea
D_{(k)}{}_{a_1\dots a_k}&=& i^{k(k-1)/2}\langle\psi,\gamma_{a_1\dots a_k}\psi\rangle.
\label{259iga}
\eea
In particular,
\be
D_{(0)}=|\psi|^2.
\ee
For spinors there is the following non-trivial Pauli identity
\cite{zhelnorovich2019}: for any spinors $\xi,\chi,\varphi$ and $\psi$ there holds
\be
\langle\xi,\psi\rangle\langle\varphi,\chi\rangle
=2^{-\nu}\Biggl\{
\langle\xi,\chi\rangle\langle\varphi,\psi\rangle
+\sum_{k=1}^{{2\nu}}\frac{1}{k!}(-1)^{k(k-1)/2}
\langle\xi,\gamma^{a_1\dots a_k}\chi\rangle
\langle\varphi,\gamma_{a_1\dots a_k}\psi\rangle
\Biggr\}.
\ee
In particular, if $\psi=\varphi$ then we get an identity
for a given spinor $\psi$ and arbitrary spinors $\xi$ and $\chi$,
\be
\la\xi,\psi\ra\la\psi,\chi\ra
=2^{-\nu}\left\{
\la\xi,\chi\ra |\psi|^2
+ \sum_{k=1}^{2\nu}\frac{1}{k!}i^{k(k-1)/2}
\la\xi,\gamma^{a_1\dots a_k}\chi\ra
D_{(k) a_1\dots a_k}
\right\}.
\ee
This identity show that a spinor is determined by the collection of tensors
$D_{(k)}$ up to a phase factor.
Furthermore,
in the case $\xi=\chi=\psi$ we have
\be
\sum_{k=1}^{{2\nu}}\frac{1}{k!}
D^{a_1\dots a_k}_{(k)}
D_{(k)}{}_{a_1\dots a_k}
=\left(2^\nu-1\right)|\psi|^4.
\ee
This is an important equation. It means that
all these forms $D_{(k)}$ are bounded by the norm square
$|\psi|^4$.

By using the identity (valid in dimension $n=4$)
\bea
4\la\xi,\psi\ra\la\chi,\varphi\ra
&=&
\la\xi,\varphi\ra\la\chi,\psi\ra
+\la\xi,\gamma^{{}a_1{}}\varphi\ra\la\chi,\gamma_{{}a_1{}}\psi\ra
-\frac{1}{2}
\la\xi,\gamma^{{}a_1{}{}a_2{}}\varphi\ra\la\chi,\gamma_{{}a_1{}{}a_2{}}\psi\ra
\nn\\
&&
-\frac{1}{6}
\la\xi,\gamma^{{}a_1{}{}a_2{}{}a_3{}}\varphi\ra\la\chi,\gamma_{{}a_1{}{}a_2{}{}a_3{}}\psi\ra
+\frac{1}{4!}
\la\xi,\gamma^{{}a_1{}\dots\mu_4}\varphi\ra\la\chi,\gamma_{{}a_1{}\dots\mu_4}\psi\ra
\nn\\
\eea
we obtain for right spinors
\bea
\la\xi_+,\gamma^{{}a_1{}{}a_2{}}\varphi_+\ra\la\chi_+,\gamma_{{}a_1{}{}a_2{}}\psi_+\ra
=-8\la\xi_+,\psi_+\ra\la\chi_+,\varphi_+\ra
+4\la\xi_+,\varphi_+\ra\la\chi_+,\psi_+\ra.
\label{2103iga}
\eea


We will be using the chiral forms
\bea
D^\pm_{(k)}{}_{a_1\dots a_k}
&=& i^{k(k-1)/2}\langle\psi_\pm,\gamma_{a_1\dots a_k}\psi_\pm\rangle.
\eea
\begin{proposition}
The chiral forms of odd degrees vanish,
\be
D^\pm_{(2k+1)}=0.
\label{262iga}
\ee
\end{proposition}
\proof
It is easy to see that the chiral forms satisfy
\bea
\langle\psi_\pm,\gamma_{a_1\dots a_{k}}\psi_\pm\rangle
&=&\langle\psi_\pm,\Gamma\gamma_{a_1\dots a_k}\Gamma\psi_\pm\rangle
=(-1)^k\langle\psi_\pm,\gamma_{a_1\dots a_k}\psi_\pm\rangle.
\label{33igax}
\eea
Therefore, they vanish for odd $k$.

The dual forms
\be
D^*_{(k)}{}_{a_1\cdots a_{{2\nu}-k}}=
\frac{1}{k!}
\varepsilon_{b_1\cdots b_{k}a_1\cdots a_{{2\nu}-k}}
D_{(k)}^{b_1\cdots b_{k}}
\ee
have the form
\be
D^*_{(k)}{}_{a_1\cdots a_{{2\nu}-k}}
=
(-1)^{k}i^{\nu-k(k-1)/2}
\langle\psi,\Gamma{}\gamma_{a_1\dots a_{{2\nu}-k}}\psi\rangle.
\ee
In particular,
\be
D^*_{(2\nu)}= \langle\psi,\Gamma\psi\rangle
=|\psi_+|^2-|\psi_-|^2.
\ee
For the chiral forms this means
\be
D^{\pm *}_{(2\nu)}= \pm |\psi_\pm|^2.
\ee
\begin{proposition}
The chiral forms of even degrees
satisfy the duality condition
\bea
D^{\pm *}_{(2k)}
&=& \pm D^{\pm}_{(2\nu-2k)}.
\eea
\end{proposition}
\proof
The dual of the chiral forms are
\bea
D^{\pm *}_{(2k)}{}_{a_1\cdots a_{2\nu-2k}}
&=&\pm
i^{\nu-k(2k-1)}
\langle\psi_\pm,{}\gamma_{a_1\dots a_{2\nu-2k}}\psi_\pm\rangle.
\eea
The statement follows.

In four dimensions, $n=4$, the only non-zero chiral forms are
$D^\pm_{(0)}$, $D^\pm_{(2)}$ and $D^\pm_{(4)}$, which is dual to $D^\pm_{(0)}$,
\be
D^\pm_{(4)abcd}=\pm\varepsilon_{abcd}D^\pm_{(0)},
\ee
with the norm
\be
|D^\pm_{(4)}|^2=\frac{1}{4!}D^\pm_{(4)abcd}D^\pm_{(4)}{}^{abcd}
=|\psi_\pm|^4
\ee
They satisfy the identity
\be
|D^\pm_{(2)}|^2=2|\psi_\pm|^4.
\ee


We consider {\it complex} {\it spinor-valued} $p$-forms;
we just call them $p$-forms without specifying that they
are spinor valued for the sake of simplicity;
this should not cause any confusion.
In other words, we consider the vector space
$S_p=\Lambda_p\otimes S$,
where $\Lambda_p$
is the space of real-valued $p$-forms.
This space has the ({ complex})
dimension
\be
\dim S_p=2^{\nu}\binom{2\nu}{p}
\ee
(recall that $n=2\nu$).
We introduce the inner product on $S_p$ by
\be
\left<\psi,\varphi\right>=\frac{1}{p!}
\bar\psi_{a_1\dots a_p}
\varphi^{a_1\dots a_p}.
\ee

%

We define the total vector space of all forms, $S_*$, by
\be
S_*=\bigoplus\limits_{p=0}^{{2\nu}}S_{p}.
\ee
The elements of the space $S_*$ are vectors containing
spin-valued forms of all degrees.
The complex dimension of this space is
\be
\dim S_*=2^{3\nu}.
\ee



There are four involutions on these spaces: $\Gamma, \zeta$,
$\tau_1$ and $\tau_2$.
First, we have the {\it chirality involution} $\Gamma$,
which decomposes the space $S_*$ into the right
and left spinors,
\be
{S}_*={S}_+\oplus{S}_-.
\ee
Second, there is the {\it $Z_2$-grading} $\zeta$ of $p$-forms
defined by
\be
\zeta_p\psi_p=(-1)^p\psi_p.
\ee
This enables one to decompose the space $S_*$ by
\be
S_*=S_{\rm even}\oplus S_{\rm odd},
\ee
where (recall that $n=2\nu$)
\be
S_{\rm even}=\bigoplus\limits_{p=0}^{\nu}S_{2p},
\qquad
S_{\rm odd}=\bigoplus\limits_{p=0}^{\nu-1} S_{2p+1},
\ee
The Hodge star $*$ acts on the spinor valued forms as usual,
\be
*_p: S_p\to S_{{2\nu}-p},
\ee
with  the properties 
\be
*_{{2\nu}-p}*_{p}=(-1)^pI_p=\zeta_p I_p,
\ee
and
\be
*(\bar\psi\wedge *\varphi) = \la\psi,\varphi\ra.
\ee
Therefore,
\be
\la\psi_{2\nu-p},*_p\varphi_p\ra
=(-1)^p\la*_{2\nu-p}\psi_{2\nu-p},\varphi_p\ra
\ee
and
\be
\la *\psi,*\varphi\ra
=\la\psi,\varphi\ra.
\ee
This enables us to define other maps
$\tau_{1}, \tau_{2}: S_p\to S_{n-p}$
by
\be
\tau_{1,p}=i^{p(p+1)+\nu}*_p,
\qquad
\tau_{2,p}=i^{p(p-1)+\nu}*_p,
\ee
(recall that $\nu=n/2$) with the property
\be
\tau_{j,2\nu-p}\tau_{j,p}=I_p.
\ee
Therefore, it can be extended to {\it duality involutions}
\be
\tau_1,\tau_2: S_*\to S_*
\ee
on the space $S_*$.
These involutions have the eigenvalues $+1$ and $-1$ and define
other decomposition of the space $S_*$ into self-dual
and anti-self-dual vectors
\be
S_*=S_{\rm sd}\oplus S_{\rm asd}.
\ee




Every spinor-valued $p$-form $\psi_{{}c_1{}\dots c_p}$
defines a collection of tensors
\be
D_{(k)}{}^{c_1\dots c_k}{}^{a_1\dots a_p}{}_{b_1\dots b_p}
= i^{k(k-1)/2}\bar\psi^{a_1\dots a_p}\gamma^{{}c_1{}\dots{}c_k{}}\psi_{b_1\dots b_p}.
\label{293iga}
\ee
These tensors are complex, in general, and
satisfy the following conjugation relation
\be
\overline{\left(D_{(k)}{}^{{}c_1{}\dots{}c_k{}}{}^{a_1\dots a_p}{}_{b_1\dots b_p}\right)}
=D_{(k)}{}^{{}c_1{}\dots{}c_k{}}{}_{b_1\dots b_p}{}^{a_1\dots a_p}{},
\ee
For $p$-forms $\psi^{b_1\dots b_p}_+$
valued in right spinors, satisfying the condition
\be
\Gamma{}\psi^{b_1\dots b_p}_+=\psi^{b_1\dots b_p}_+,
\ee
the tensors of odd degree vanish,
\be
D^+_{(2k+1)}{}^{{}c_1{}\dots c_{2k+1}}{}^{b_1\dots b_p}{}_{a_1\dots a_p}=0.
\label{294iga}
\ee


We restrict ourselves to four dimensions, $n=2\nu=4$,
and to the right spinor subspace
\be
S_*^+=\Ker P_-=\{\psi\in S_*\;|\; \Gamma\psi=\psi\}.
\ee
Recall that for the right spinors the following tensors vanish,
(by using (\ref{294iga})),
\be
\bar\psi_+^{ {}a_1{}\dots a_{p}}\gamma_c\psi^+_{b_1\dots b_{p}}
=\bar\psi_+^{ {}a_1{}\dots a_{p}}\gamma_{cdf}\psi^+_{b_1\dots b_{p}}=0.
\ee
Therefore, the only non-zero tensors constructed from the spinor-valued forms
are the tensors $D_{(0)}$ and $D_{(2)}$
defined in (\ref{2106iga}) and (\ref{2107iga})
(and, of course, $D_{(4)}$ which is dual to $D_{(0)}$).
For the right spinor valued $p$-forms the nonzero tensors in four dimensions are,
$p=0,1,2,3,4$,
\bea
D^+_{(0)}{}^{a_1\dots a_p}{}_{b_1\dots b_p}
&=& \bar\psi_+^{a_1\dots a_p}\psi_{+ b_1\dots b_p},
\label{2106iga}
\\
D^+_{(2)}{}^{cd}{}^{a_1\dots a_p}{}_{b_1\dots b_p}
&=& i\bar\psi_+^{a_1\dots a_p}\gamma_+^{cd}\psi_{+ b_1\dots b_p},
\label{2107iga}
\\
D^+_{(4)}{}^{rstu \;a_1\dots a_p}{}_{b_1\dots b_p}
&=&E^{rstu}
D^+_{(0)}{}^{a_1\dots a_p}{}_{b_1\dots b_p}.
\eea
This defines a real 2-form for $p=0$
\be
D^+_{(2)}{}^{ab}
= i\bar\psi_+\gamma_+^{ab}\psi_+,
\label{2112iga}
\ee
and a complex tensor for $p=1$
\bea
D^+_{(2)}{}^{cd}{}^{ab}
&=& i\bar\psi_+^{a}\gamma_+^{cd}\psi_+^{b}
\label{2113iga}
\eea


By using eq. (\ref{2103iga}) we obtain the identity
\bea
D^+_{(2)}{}^{rs}{}^{a_1\dots a_p}{}_{b_1\dots b_p}
D^+_{(2)}{}_{rs}{}^{c_1\dots c_p}{}_{d_1\dots d_p}
=
8D^+_{(0)}{}^{}{}^{a_1\dots a_p}{}_{d_1\dots d_p}
D^+_{(0)}{}{}^{c_1\dots c_p}{}_{b_1\dots b_p}
\nn\\
-4D^+_{(0)}{}^{}{}^{a_1\dots a_p}{}_{b_1\dots b_p}
D^{+}_{(0)}{}{}^{c_1\dots c_p}{}_{d_1\dots d_p}.
\label{2122iga}
\eea
For $p=0$ this gives, in particular, the norm of the tensor
$D^+_{(2)\mu\nu}$,
\bea
D^+_{(2)}{}^{ab}D^+_{(2)}{}_{ab}=4|\psi_+|^4.
\label{450uni}
\eea

In four dimensions the spin group is the product
\be
\Spin(4)=SU(2)\times SU(2).
\ee
Therefore, the irreducible representations $V_{(i,j)}$
 of the spin group
$\Spin(4)$ are labeled by two half-integers $(i,j)$
of dimension 
\be
\dim V_{(i,j)}=(2i+1)(2j+1),
\ee
The fundamental representations are the spinor
representation and the cospinor representations
\be
S^+_0=\left(\frac{1}{2},0\right),
\qquad
S^-_0=\left(0,\frac{1}{2}\right),
\ee
 corresponding to the right and left Weyl
spinors; of course, the dimensions of these spaces are
\be
\dim S^\pm_0=2.
\ee

 The one forms transform according to the representation 
 $\Lambda_1=\left(\frac{1}{2}, \frac{1}{2}\right)$. Therefore, right and the left spinor valued one forms decompose by
\bea
S_1^+ &=&\Lambda_1\otimes S^+_0
=\left(1,\frac{1}{2}\right)
\oplus
\left(0,\frac{1}{2}\right),
\qquad
\\
S_1^-&=&\Lambda_1\otimes S^-_0
=\left(\frac{1}{2},1\right)
\oplus\left(\frac{1}{2},0\right).
\eea
The dimensions of these representations are
\be
\dim S_1^\pm = 2\cdot 3+2=8.
\ee
Here $\left(1,\frac{1}{2}\right)$ and $\left(\frac{1}{2},1\right)$ are gamma traceless
right and left  one forms.
That is, right (left) spinor valued one forms
are decomposed according to
\be
\psi^\pm_\mu=\varphi^\pm_\mu+\gamma_\mu \xi^\mp,
\ee
where
$\varphi^\pm_\mu$
satisfy the conditions
\be
\Gamma\varphi^\pm_\mu=\pm\varphi^\pm_\mu,
\qquad
\gamma^\mu\varphi^\pm_\mu=0,
\ee 
 and
\be
\xi^\mp=\frac{1}{4}\gamma^\mu\psi_\mu^\pm.
\ee

The self-dual two forms 
 and anti-self-dual two-forms 
transform by the representations
\be
\Lambda^{+}_{2}=(1,0),
\qquad
\Lambda^{-}_2=(0,1).
\ee
 Therefore,
the right and left spinor valued two forms decompose as follows:
\bea
S_{2}^{+, sd}&=&\Lambda^{+}_{2}\otimes S^+_0
=\left(\frac{1}{2},0\right)\oplus 
\left(\frac{3}{2},0\right),
\\
S_{2}^{+, asd}&=&\Lambda^{-}_{2}\otimes S^+_0
=\left(\frac{1}{2},1\right),
\\
S_{2}^{-, sd}&=&\Lambda^{+}_{2}\otimes S^-_0
=\left(1,\frac{1}{2}\right),
\\
S_{2}^{-, asd}&=&\Lambda^{-}_{2}\otimes S^-_0
=\left(0,\frac{1}{2}\right)\oplus 
\left(0,\frac{3}{2}\right).
\eea

Therefore,
\bea
S_{2}^{+}&=&
\left(\frac{3}{2},0\right)\oplus
\left(\frac{1}{2},1\right)\oplus
\left(\frac{1}{2},0\right),
\\
S_{2}^{-}&=&
\left(0,\frac{3}{2}\right)
\oplus\left(1,\frac{1}{2}\right)\oplus
\left(0,\frac{1}{2}\right).
\eea
The dimensions of these representations are
\be
\dim S_2^\pm = 4 + 2\cdot 3+ 2 = 12.
\ee

Here
$\left(\frac{3}{2},0\right)$
corresponds to the self-dual right 
spinor valued gamma traceless 
two forms satisfying the conditions
\bea
\Gamma\psi^+_{\mu\nu}=(*\psi^+)_{\mu\nu}
=\psi^+_{\mu\nu},
\qquad
\gamma^\mu\psi^+_{\mu\nu}=0,
\eea
and
$\left(0,\frac{3}{2}\right)$ corresponds
to the anti-self-dual left  gamma traceless
spinor valued two forms
satisfying the conditions
\bea
\Gamma\psi^-_{\mu\nu}=(*\psi^-)_{\mu\nu}=-\psi^-_{\mu\nu},
\qquad
\gamma^\mu\psi^-_{\mu\nu}=0.
\eea

Every right (left) spinor valued two form has the form
\be
\psi^\pm_{\mu\nu}=\varphi^\pm_{\mu\nu}
+2\gamma_{[\mu}\xi^\mp_{\nu]}
+\gamma_{\mu\nu}\eta^\pm,
\ee
where
\be
\gamma^\mu\varphi^\pm_{\mu\nu}=\gamma^\mu\xi^\pm_\mu=0,
\ee
and
\bea
\eta^\pm &=& -\frac{1}{12}\gamma^{\mu\nu}\psi^\pm_{\mu\nu},
\\
\xi_\nu^\mp &=& \frac{1}{2}\gamma^\mu\psi^\pm_{\mu\nu}
+\frac{1}{8}\gamma_\nu\gamma^{ab}\psi^\pm_{ab}.
\eea


\section{Algebra of Spinor Valued Forms}
\setcounter{equation}0
\label{sec:algebra}


Let $\theta^a$, $a=1,\dots, 2\nu$, be the Grassmann variables satisfying
the anticommuting relations
\be
\theta^a\theta^b=-\theta^b\theta^a.
\ee
We will identify the elements of the total space of spinor valued forms
$S_*$ with functions of anti-commuting variables $\theta^a$,
that is, every $\psi\in S_*$ has the form
\be
\psi(x,\theta)=\sum_{p=0}^{2\nu} \frac{1}{p!}\psi_{a_1\dots a_p}(x)
\theta^{a_1}\cdots\theta^{a_p}.
\ee
We can define the inner product in this vector space by
\be
\la\psi,\varphi\ra=\sum_{p=0}^{2\nu}\la\psi_p,\varphi_p\ra,
\ee
where $\psi_p$ and $\varphi_p$ are the $p$-components of the general
vectors $\psi$ and $\varphi$.
This inner product can also be defined in terms of the Berezin integral
over the Grassmann variables.
The Hodge star $*$ as well as the involutions $\Gamma$, $\zeta$, $\tau_1$ 
and $\tau_2$ act on 
the space $S_*$ accordingly.


Further, we define the
(left)  derivative by
\be
\frac{d}{dt}\psi(\theta+t h)\Big|_{t=0}
=h^a\frac{\partial}{\partial \theta^a}\psi(\theta);
\ee
then
\be
\frac{\partial}{\partial \theta^a}\theta^b=\delta^b_a.
\ee
and, more generally,
\bea
\frac{\partial}{\partial \theta^a}\theta^{b_1}\cdots\theta^{b_p}
&=&p\delta^{[b_1}_a\theta^{b_2}\cdots\theta^{b_p]}
\eea
The derivatives are anti-commuting operators
\be
\frac{\partial}{\partial \theta^a}\frac{\partial}{\partial \theta^b}
=-\frac{\partial}{\partial \theta^b}\frac{\partial}{\partial \theta^a}
\ee
The derivative (anti)-commutes with the multiplication
operator as follows
\be
\frac{\partial}{\partial \theta^a}\theta^b
+\theta^b\frac{\partial}{\partial \theta^a}
=\delta^b_a.
\label{243igax}
\ee

Notice that the operators $X_a\theta^a$ and $X^a\frac{\partial}{\partial \theta^a}$
act on $\psi$ by
\bea
X_a\theta^a\psi &=& \sum_{p=0}^{n-1} \frac{(p+1)}{(p+1)!}X_{[a}\psi_{b_1\dots b_{p}]}
\theta^a\theta^{b_1}\cdots\theta^{b_{p}},
\\
X^a\frac{\partial}{\partial \theta^a}\psi
&=& \sum_{p=0}^{n-1} \frac{1}{p!}X^a\psi_{a b_1\dots b_{p}}
\theta^{b_1}\cdots\theta^{b_{p}},
\eea
and define the exterior product 
and the interior product with the one-form $X$.



We define the operators $\alpha,\beta, h, \varkappa: S_*\to S_*$
on the total space $S_*$ by
\be
\alpha=\gamma_a \theta^a,
\qquad
\beta=\gamma^a\frac{\partial}{\partial \theta^a},
\qquad
h=\theta^a\frac{\partial}{\partial \theta^a},
\label{289iga}
\ee
and 
\be
\varkappa=2\nu I- 2h.
\ee
The restrictions of these operators to the spaces $S_p$
define the maps
\bea
\alpha_p : {S}_p\to {S}_{p+1},
\qquad
\beta_p : {S}_p\to {S}_{p-1},
\eea
by: for any $\psi\in S_p$,
\bea
(\alpha_p\psi)_{a_1\dots a_{p+1}}&=&
(p+1)\gamma_{[a_1}\psi_{a_2\dots a_{p+1}]},
\\
(\beta_p\psi)_{a_1\dots a_{p-1}}&=&
\gamma^b\psi_{b a_1\dots a_{p-1}},
\eea
Obviously,
\be
\alpha_{2\nu}=0,
\qquad
\beta_0=0,
\qquad
\varkappa_{\nu}=0.
\ee

The powers of the operators $\alpha$ and $\beta$
\bea
\alpha_p^k : {S}_p\to {S}_{p+k}, \qquad
\beta_p^k : {S}_p\to {S}_{p-k}
\eea
have the form
\bea
(\alpha_p^k\psi)_{a_1\dots a_{p+k}}&=&
\frac{(p+k)!}{p!}\gamma_{[a_1\dots a_k}\psi_{a_{k+1}\dots a_{p+k}]},
\\
(\beta_p^k\psi)_{a_1\dots a_{p-k}}&=&
(-1)^{k(k-1)/2}\gamma^{b_1\dots b_k}\psi_{b_1\dots b_k a_1\dots a_{p-k}}.
\eea
Obviously, these operators are nilpotent since they have only finitely many
non-trivial powers
\be
\alpha_p, \ \alpha_p^2, \dots, (\alpha_p)^{2\nu-p},
\qquad
\beta_p, \ \beta_p^2, \dots, (\beta_p)^p\,.
\ee
that is,
\be
\alpha^k=\beta^k=0, \qquad \mbox{for } k>2\nu.
\ee


\begin{proposition}
	The operators $\alpha_p$ and $\beta_p$ satisfy
	the commutation relation
	\be
	\beta_{p+1}\alpha_p-\alpha_{p-1}\beta_p=2(\nu-p)I_p.
	\label{285fg}
	\ee
	Therefore, for $p=\nu$ they commute, that is,
	$
	\alpha_{\nu-1}\beta_\nu=\beta_{\nu+1}\alpha_\nu.
	$
\end{proposition}
\proof
We have for any $\psi\in S_p$
\bea
(\alpha_{p-1}\beta_p\psi)_{a_1\dots a_{p}}
&=&
p\gamma_{[a_1}\gamma^b\psi_{|b|a_2\dots a_{p}]},
\\
(\beta_{p+1}\alpha_p\psi)_{a_1\dots a_{p}}
&=&
(p+1)\gamma^b\gamma_{[b}\psi_{a_1\dots a_{p}]}.
\eea
Therefore,
\be
(p+1)\gamma^b\gamma_{[b}\psi_{a_1\dots a_{p}]}
=(2\nu-2p)\psi_{a_1\dots a_{p}}
+p\gamma_{[a_1}\gamma_{|b|}\psi^b{}_{a_2\dots a_{p}]},
\ee
which gives (\ref{285fg}).

\begin{proposition}
	The operators $\alpha$, $\beta$ and $\varkappa$ form the
	Lie algebra $\slfrak(2,\RR)$
	\bea
	[\beta,\alpha]&=&\varkappa,
	\label{247igax}\\
	{}[\varkappa,\alpha] &=&-2\alpha,
	\label{248igax}\\
	{}[\varkappa,\beta] &=& 2\beta.
	\label{249igax}\eea
	The Casimir operator of the Lie algebra,
	\bea
C&=&\frac{1}{2}(\alpha\beta+\beta\alpha)+\frac{1}{4}\varkappa^2
	\label{299iga}\eea
	commutes with all elements of the algebra
	\be
	[C,\alpha]=[C,\beta]=[C,\varkappa]=0.
	\ee
\end{proposition}
\proof
By using (\ref{12igax}) and (\ref{243igax}) we obtain (\ref{247igax}).
Further, by computing the commutators
\bea
[h,\theta^b] &=&  \theta^b,
\label{292iga}\\
{}\left[h, \frac{\partial}{\partial \theta^b}\right]
&=& -\frac{\partial}{\partial \theta^b}.
\label{293igax}
\eea
we get (\ref{248igax}) and (\ref{249igax}).


That is why, the operator $\alpha$ is  the lowering operator,
the operator $\beta$ is the raising operator and the operator
$\varkappa$ is the Cartan operator of the algebra
$\slfrak(2,\RR)$.
The operator $h$ acts on homogeneous functions of degree $p$
\be
\psi_p=\frac{1}{p!}\psi_{a_1\dots a_p}(x)
\theta^{a_1}\cdots\theta^{a_p},
\ee
by
\be
h\psi_p=p\psi_p.
\ee
Therefore, the Cartan operator acts on this functions
diagonally
\be
\varkappa\psi_p=(2\nu-2p)\psi_p.
\ee

\begin{proposition}
	The Cartan operator is Hermitian 
	\be
	\varkappa^*=\varkappa.
	\ee
and	anti-commutes with the Hodge star
	\be
	*\varkappa=-\varkappa*.
	\ee
\end{proposition}
\proof
This follows from the relation $\varkappa=2\nu-2h$ and
\be
*h=(2\nu-h)*.
\ee

\begin{proposition}
The operators $\beta$ and $\alpha$ are the adjoints of each other,
\be
	\beta=\alpha^*
\ee
and are intertwined by the Hodge star
	\be
	*_{p-1}\beta_p = (-1)^{p+1}\alpha_{n-p} *_p.
	\ee
\end{proposition}
\proof
For any $\varphi\in S_p, \psi\in S_{p+1}$,
we have
\be
\la\psi,\alpha_p\varphi\ra
=\la\beta_{p+1}\psi,\varphi\ra.
\ee
Therefore, $\beta=\alpha^*$.
The duality relation follows from the equation
\be
\gamma_{d} \psi^{d c_1 \dots c_{p-1}}
=\frac{1}{p!(n-p)!}
\varepsilon^{b_1b_2\dots b_{2\nu-p+1}c_1\dots c_{p-1}}
\gamma_{b_1}
\varepsilon_{a_1\dots a_p b_2\dots b_{2\nu-p+1}}
\psi^{a_1\dots a_p}.
\ee


We define the adjoint operator $\ad_X$ on the algebra as usual
\be
\ad_XY=[X,Y].
\ee
Then we have the identities
\bea
\ad_\alpha^2\varkappa =
\ad_\beta^2 \varkappa  = 0.
\eea
and for any $k\ge 1$ there hold
\bea
\ad_\varkappa^k\alpha  &=& (-2)^k\alpha,
\\
\ad_\varkappa^k \beta &=& 2^k\beta.
\eea
Furthermore,
\bea
\ad_\alpha^2\beta  &=& -2\alpha
\\
\ad_\beta^2 \alpha &=& -2\beta,
\eea
and
for any $k\ge 3$
\bea
\ad_\alpha^k\beta  =
\ad_\beta^k \alpha = 0,
\eea


\begin{corollary}
	There hold
	\bea
	\exp(t\alpha)\beta\exp(-t\alpha) &=& \beta -t\varkappa -t^2\alpha,
	\\{}
	\exp(t\beta)\alpha\exp(-t\beta) &=& \alpha +t\varkappa -t^2\beta,
	\\
	\exp(t\alpha)\varkappa\exp(-t\alpha) &=& \varkappa + 2t\alpha,
	\\
	\exp(t\beta)\varkappa\exp(-t\beta) &=& \varkappa - 2t\beta,
	\\
	\exp(t\varkappa)\beta\exp(-t\varkappa) &=& e^{2t}\beta,
	\\{}
	\exp(t\varkappa)\alpha\exp(-t\varkappa) &=& e^{-2t}\alpha.
	\eea
	
\end{corollary}
\proof
This is proved by using the series
\be
e^{tX} Y e^{-tX}=\sum_{k=0}^\infty\frac{t^k}{k!}\ad^k_XY.
\ee

\begin{proposition}
	There hold
	\bea
	[\beta,\alpha^k] &=& \left(\varkappa+k-1\right)k\alpha^{k-1},
	\label{2114iga}\\{}
	[\beta^k,\alpha] &=& \left(\varkappa
	-k+1\right)k\beta^{k-1},
	\label{322fg}\\{}
	[\alpha^k,\varkappa] &=& 2k\alpha^k,
	\label{259igax}
	\\{}
	[\beta^k,\varkappa] &=& -2k\beta^k,
	\label{275iga}\\{}
	\alpha\varkappa^k &=& (\varkappa+2)^k\alpha,
	\\{}
	\beta\varkappa^k &=& (\varkappa-2)^k\beta.
	\eea
\end{proposition}
\proof
This is easily proved by using the equations above.

We define
the operators $N_+, N_-: S_*\to S_*$ by
\be
N_+=\alpha\beta,
\qquad
N_-=\beta\alpha,
\ee
and the operator
\bea
H&=& N_++N_-=\alpha\beta+\beta\alpha.
\eea
The restrictions of these operators
$N_p^+, N_p^-: S_p\to S_p$
have the form
\be
N_p^+=\alpha_{p-1}\beta_p,
\qquad
N_p^-=\beta_{p+1}\alpha_p.
\ee
They are related by the intertwining relation
\be
N_-*=*N_+,
\ee
or
\be
*_{p}N_p^+=N_{n-p}^-*_{p}.
\ee
Notice that in the middle dimension they are equal
\be
N^+_\nu=N^-_\nu.
\ee

\begin{proposition}
	The operators $N_+$ and $N_-$ satisfy the commutation
	commutation relations
	\bea
	[N_+,\alpha] &=& (\varkappa+2)\alpha,
	\\{}
	[N_+,\beta] &=& -\varkappa\beta.
	\\{}
	[N_-,\alpha] &=& \varkappa\alpha,
	\\{}
	[N_-,\beta] &=& -(\varkappa-2)\beta,
	\\{}
	[N_+,N_-]&=& [N_+,\varkappa]=[N_-,\varkappa]=0.
	\eea
	More generally, for any positive integer $k\ge 1$
	\bea
	[N_+, \alpha^k] &=& (\varkappa+k+1)k\alpha^k,
	\label{2116iga}
	\\{}
	[N_+,\beta^k] &=& -(\varkappa-k+1)k\beta^k,
	\\{}
	[N_-, \alpha^k] &=& (\varkappa+k-1)k\alpha^k,
	\\{}
	[N_-,\beta^k] &=& -(\varkappa-k-1)k\beta^k.
	\eea
\end{proposition}
\proof
This follows from above.


\begin{proposition}
	The operators
	$
	N_+
	$,
	 $
	N_-$ and $H$
	are Hermitian and non-negative. That is, for any $\psi$,
	\be
	\la\psi, N_+\psi\ra\ge 0,
	\qquad
	\la\psi, N_-\psi\ra\ge 0,
	\qquad
	\la\psi, H\psi\ra\ge 0.
	\ee
\end{proposition}
\proof
This follows from the equations
\be
\la\psi, N_+\psi\ra=|\beta\psi|^2,
\qquad
\la\psi, N_-\psi\ra=|\alpha\psi|^2,
\qquad
\la\psi, H\psi\ra=|\alpha\psi|^2+|\beta\psi|^2.
\ee

\begin{corollary}
	The kernels of the operators $N_+, N_-$ and $H$ are determined by the kernels
	of the operators $\alpha$ and $\beta$,
	\be
	\Ker N_+=\Ker\beta,
	\qquad
	\Ker N_-=\Ker\alpha,
	\qquad
	\Ker H=\Ker\beta\cap \Ker \alpha.
	\ee
\end{corollary}


Notice that
\be
\Ker\alpha_{2\nu}=S_{2\nu},
\qquad
\Ker\beta_0=S_0,
\qquad
\Ker\varkappa_\nu=S_\nu.
\ee

\begin{proposition}
For $p=0$ the (complex) dimension of the kernel of the operator $\beta_0$
is 
\be
\dim\Ker\beta_0=2^{\nu}.
\ee
For $p=1,\dots,\nu$,
\bea
\dim\Ker\beta_p
=2^{\nu}
{2\nu\choose p}\frac{2\nu-2p+1}{2\nu-p+1}.
\eea
\end{proposition}
\proof
For $p=0$ we have $\Ker\beta_0=S_0$. This gives the dimension of the kernel.

For $p=1,\dots,\nu$, the map $\beta_p: S_p\to S_{p-1}$ is surjective. 
Therefore, the dimension of the image is exactly the dimension of the target space,
\bea
\dim\Ker\beta_p = 2^{\nu}\left\{
{2\nu\choose p} -{2\nu\choose p-1}
\right\}.
\eea


\begin{proposition}{}
	\label{prop5iga}
	\benum
	\item
	If $p>\nu$ then
	$
	\Ker\beta_p=\Ker N^+_p=\{0\},
	$
	\item
	If $p<\nu$ then
	$
	\Ker\alpha_p=\Ker N^-_p=\{0\},
	$
	\item
	If $p=\nu$ then
	\be
	\Ker\alpha_\nu=\Ker\beta_\nu.
	\ee
	\eenum
\end{proposition}
\proof
This follows from the equation, for any form $\varphi\in S_p$,
\be
|\alpha_p\varphi|^2=|\beta_p\varphi|^2+2(\nu-p)|\varphi|^2.
\label{2127iga}\ee

\begin{corollary}{}
The kernel of the operator $\beta$ has the form
	\be
	\Ker\beta
	=S_0\oplus\Ker\beta_1\oplus\cdots\oplus\Ker\beta_{\nu}.
	\ee
\end{corollary}

\begin{proposition}
The kernel of the operator
$H_p$ is trivial except in the middle
dimension, that is, if $p\ne\nu$ then
\be
\Ker H_p=\{0\}.
\ee
\end{proposition}
\proof
Since $\Ker\beta_p=\{0\}$ for $p>\nu$ and
$\Ker \alpha_p=\{0\}$ for $p<\nu$  then
$
\Ker H_p=\Ker\beta_p\cap \Ker \alpha_p
=\{0\}$ for $p\ne\nu$.
Alternatively, if $\psi\in\Ker H_p$
then
\be
\beta_p\psi=\alpha_p\psi=0.
\ee
By using (\ref{285fg}) this means
\be
(\beta_{p+1}\alpha_p-\alpha_{p-1}\beta_p)\psi
=\varkappa_p\psi=2(\nu -p)\psi=0.
\ee
Therefore, $\psi=0$ for any $p\ne \nu$.

The following subspace plays a special role
\be
\cH=\Ker H_\nu=\Ker N^-_\nu=\Ker N^+_\nu
=\Ker\alpha_\nu=\Ker\beta_\nu \subseteq S_\nu.
\ee
It can be canonically decomposed into the subspaces of self-dual
and antiself-dual forms,
\be
\cH=\cH_{\rm sd}\oplus\cH_{\rm asd}.
\ee

\begin{proposition}
	The spaces
	$\Imag \alpha_{p-1}$
	and $\Ker\beta_{p}$
	are orthogonal.
\end{proposition}
\proof
We have
for any $\psi\in \Ker\beta_p$, $\varphi\in \Imag\alpha_{p-1}$
\be
\la\psi,\varphi\ra=\la\psi,\alpha_{p-1}\xi\ra=\la \beta_{p}\psi,\xi\ra=0.
\ee
The statement follows.

\begin{proposition}
There are the orthogonal decompositions
\be
S_p=\Ker\beta_p\oplus\Imag\alpha_{p-1}=\Ker\alpha_p\oplus\Imag\beta_{p+1}.
\ee
In particular,
\benum
\item
if $p<\nu$, then $\beta_{p+1}$ is surjective, $\Imag\beta_{p+1}=S_p$, and
$S_p=\Ker\beta_p\oplus\Imag\alpha_{p-1}$;
\item
if $p>\nu$, then $\alpha_{p-1}$ is surjective, $\Imag\alpha_{p-1}=S_p$, and
$S_p=\Ker\alpha_p\oplus\Imag\beta_{p+1}$;
\item
if $p=\nu$, then
$S_\nu=\cH\oplus\Imag\alpha_{\nu-1}=\cH\oplus\Imag\beta_{\nu+1}$,
with $\Imag\alpha_{\nu-1}=\Imag\beta_{\nu+1}=\cH^\perp$.
\eenum
Consequently,
\be
S_*=\cH\oplus\Imag H,
\qquad
\Imag H=\Imag\alpha+\Imag\beta,
\ee
where the last sum is not direct.
\end{proposition}
\proof
Since $\beta_p=\alpha_{p-1}^*$, we have
\be
\Imag\alpha_{p-1}=(\Ker\beta_p)^\perp,
\qquad
\Imag\beta_{p+1}=(\Ker\alpha_p)^\perp
\ee in $S_p$, which gives the two
orthogonal decompositions. By (\ref{2127iga}), $\alpha_p$ is injective for
$p<\nu$ and $\beta_p$ is injective for $p>\nu$; hence their adjoints
$\beta_{p+1}$ and $\alpha_{p-1}$ are surjective. For $p=\nu$ we have
\be
\Ker\alpha_\nu=\Ker\beta_\nu=\cH.
\ee Finally,
\be
\Imag H\subseteq\Imag\alpha+\Imag\beta
\ee since 
$H=\alpha\beta+\beta\alpha$,
and $\Imag\alpha=(\Ker\beta)^\perp\subseteq(\Ker H)^\perp=\Imag H$,
$\Imag\beta=(\Ker\alpha)^\perp\subseteq\Imag H$.


\section{Orthonormal Basis of Spinor Valued Forms}
\setcounter{equation}0
\label{sec:basis}

We construct a basis in the vector space $S_*$.
The space $S_*$ is reducible; it decomposes into irreducible 
representations of the Lie algebra $\slfrak(2,\RR)$
constructed as follows.
We define the subspaces $S_{(j,m)}\subseteq S_{\nu-m} $ by
\be
S_{(j,m)}=\alpha^{j-m}\Ker\beta_{\nu-j},
\ee
where $j=0,1,\dots,\nu$ and $m=-j,\dots,j$.

Let $\xi\in\Ker\beta$ be a nonzero vector determined by
\be
\beta\xi=0.
\ee 
This vector has all degrees $p=0,1,\dots,\nu$,
\be
\xi=\sum_{p=0}^{\nu} \frac{1}{p!}\xi_{a_1\dots a_p}
\theta^{a_1}\cdots\theta^{a_p}.
\ee

Each $p$ generates an irreducible representation.
Then
\be
\varkappa\xi_{\nu-j}=2j\xi_{\nu-j},
\ee
where
\be
j=\nu-p.
\ee
Let $\varphi_{j,m}\in S_{(j,m)}$  be the vectors
defined by
\be
\varphi_{j,m}=\alpha^{j-m}\xi_{\nu-j},
\ee
with $m=-j,\; -j+1,\;\dots,\;j-1,\;j$.
The vector
\be
\varphi_{j,j}=\xi_{\nu-j}.
\ee
is called the highest weight vector
of the irreducible representation $j$.

\begin{proposition}
	\begin{enumerate}
		\item
		The vectors $\varphi_{j,m}$ are eigenvectors of the Cartan operator $\varkappa$, 
		\be
		\varkappa\varphi_{j,m}=2m\varphi_{j,m}.
		\ee
		\item
		The vectors $\varphi_{j,m}$ are mutually orthogonal, that is, for 
any $j$ and any $m\ne m'$,
		\be
		\la\varphi_{j,m},\varphi_{j,m'}\ra=0.
		\ee
	\end{enumerate}
\end{proposition}
\proof
\begin{enumerate}
	\item
	By using (\ref{259igax}) we get
	\be
	\varkappa\varphi_{j,m}=\varkappa\alpha^{j-m}\xi_{\nu-j}
	=\alpha^{j-m}\left[\varkappa-2(j-m)\right]\xi_{\nu-j}
	=2m\alpha^{j-m}\xi=2m\varphi_{j,m}.
	\ee
	
	\item
	The orthogonality follows from the fact that they are eigenvectors of a Hermitian
	operator.
	
	\item
	We compute the norm of $\varphi_{j,m}$; we have
	\be
	|\varphi_{j,m}|^2=\la\alpha^{j-m}\xi_{\nu-j},\alpha^{j-m}\xi_{\nu-j}\ra
	=\la\xi,\beta^{j-m}\alpha^{j-m}\xi_{\nu-j}\ra.
	\ee
	
	\item
	We need to compute
	$
	\beta^{k}\alpha^{k}\xi_{\nu-j}.
	$
	By using  (\ref{2114iga}) we get 
	\be
	\beta^{k}\alpha^{k} = \beta^{k-1}\alpha^{k}\beta
	+k\beta^{k-1}(\varkappa+ k-1) \alpha^{k-1}
	\ee
	and, further, by using (\ref{275iga})
	we obtain
	\be
	\beta^{k}\alpha^{k} = \beta^{k-1}\alpha^{k}\beta
	+k(\varkappa-k+1) \beta^{k-1}\alpha^{k-1}
	\ee
	we get the recursion
	\be
	\beta^{k}\alpha^{k}\xi_{\nu-j}=k(\varkappa-k+1) \beta^{k-1}\alpha^{k-1}\xi_{\nu-j}.
	\ee
	
	\item
	This recursion can be solved explicitly
	\bea
	\beta^{k}\alpha^{k}\xi_{\nu-j} &=& k!(\varkappa-k+1)\cdots (\varkappa-1)\varkappa
	\xi_{\nu-j}
	\nn\\
	&=&
	k!2j(2j-1)\cdots (2j-k+1)\xi_{\nu-j},
	\eea
	which can be written in the form
	\bea
	\beta^{k}\alpha^{k}\xi_{\nu-j} &=&
	k!\frac{(2j)!}{(2j-k)!}\xi_{\nu-j},
	\eea
	
	\item
	This gives
	\be
	|\varphi_{j,m}|^2=\frac{(j-m)!(2j)!}{(j+m)!}|\xi_{\nu-j}|^2
	\ee

\end{enumerate}

Therefore, the orthonormal basis
is
\be
\psi_{j,m}=\sqrt{\frac{(j+m)!}{(j-m)!(2j)!}}\alpha^{j-m}
\frac{\xi_{\nu-j}}{|\xi_{\nu-j}|}.
\label{4117iga}
\ee

\begin{theorem}
	\begin{enumerate}
		\item
		The vectors $\psi_{j,m}$ are the eigenvectors of the operators
		$N_\pm$ and $H$,
		\bea
		N_+\psi_{j,m} &=& (j-m)(j+m+1)\psi_{j,m},
		\\
		N_-\psi_{j,m} &=& (j+m)(j-m+1)\psi_{j,m},
		\\
		H\psi_{j,m} &=& \left[2j(j+1)-2m^2\right]\psi_{j,m}.
		\label{292igax}
		\eea
		
		\item
		The vectors $\psi_{j,m}$ form an irreducible representation of the
		algebra $\slfrak(2,\RR)$ labeled by $j$
		of dimension $2j+1$ with the Casimir operator
		\be
		C\psi_{j,m}=j(j+1)\psi_{j,m}.
		\label{293igaz}
		\ee
		\item
		The operators $\alpha$ and $\beta$ act on the vectors $\psi_m$ by
		\bea
		\alpha\psi_{j,m}=\sqrt{(j+m)(j-m+1)}\psi_{j,m-1},
		\\
		\beta\psi_{j,m}=\sqrt{(j-m)(j+m+1)}\psi_{j,m+1},
		\eea
		\item
		The Hodge star acts on the vectors $\psi_m$ by
		\be
		*S_{(j,m)}=S_{(j,-m)}.
		\ee
	\end{enumerate}
\end{theorem}
\proof
\begin{enumerate}
	\item
	By using eq. (\ref{2114iga}) we get
	\bea
	\beta\varphi_{j,m} &=& (j-m)(m+j+1)\varphi_{j,m+1},
	\eea
	Now by using eq. (\ref{248igax}) we get
	\be
	N_+\varphi_m=(j-m)(m+j+1)\varphi_m
	\ee
	Now, by using $N_-=N_++\varkappa$ we get
	\be
	N_-\varphi_m=(j+m)(j-m+1)\varphi_m;
	\ee
	the eqs. (\ref{292igax}) follow.
	
	\item
	The dimension of this representation is equal to the number of 
	possible values of the
	integer $m=-j,-j+1,\dots, j-1,j$.
	
	\item
	The action of the operators $\alpha$ and $\beta$ on the vectors
	$\psi_{j,m}$ follows from the normalization.
	
	\item
	The duality relation follows from the anticommutativity of
	the Cartan operator $\varkappa$ with the Hodge star $*$.
	
\end{enumerate}


Since the vectors $\psi_{j,m}$ are the eigenvectors of a Hermitian
Casimir operator $C$ then they are orthogonal for $j\ne j'$.
The collection of vectors $\psi_{j,m}$ with
$j=0,\dots,\nu$ and $m=-j,\dots, j$ forms an orthonormal basis
in the space $S_*$, that is,
\be
S_*=\Span\left\{\psi_{j,m}\;|\; j=0,\dots,\nu; \quad m=-j,\dots, j
\right\}
\ee
In other words, we have the fundamental decomposition theorem.

\begin{theorem}
\be
S_*=\bigoplus_{j=0}^\nu\bigoplus_{m=-j}^j\alpha^{j-m}\Ker\beta_{\nu-j}.
\ee
\end{theorem}

As a simple check we verify the dimension of this space
\bea
\dim S_* &=& \sum_{j=0}^\nu(2j+1)\dim \Ker\beta_{\nu-j}
\\
&=&
2^{\nu}(2\nu+1)
+2^{\nu}\sum_{j=0}^{\nu-1}(2j+1)\left\{
{2\nu\choose \nu-j}
-{2\nu\choose \nu-j-1}
\right\}
\nn\\
\nn\\
&=&
2^{\nu}\sum_{j=0}^{2\nu}{2\nu\choose j}
=2^{3\nu}.
\eea

In four dimensions, $\nu=2$, there are three irreducible representations,
$j=0,1,2$,
\bea
V_0 &=& \Ker\beta_{2},
\\
V_1 &=& \alpha^2\Ker\beta_{1}\oplus \alpha^1\Ker\beta_{1}\oplus \Ker\beta_{1},
\\ 
V_2 &=& \alpha^4\Ker\beta_{0}\oplus \alpha^3\Ker\beta_{0}\oplus
\alpha^2\Ker\beta_{0}\oplus  \alpha^1\Ker\beta_{0}\oplus 
\Ker\beta_{0}.
\eea
Notice that these subspaces correspond to the following
decomposition in irreducible representations of the
group $\Spin(4)$,
\bea
\Ker\beta_0 &=&\left(\frac{1}{2},0\right)\oplus \left(0,\frac{1}{2}\right),
\\
\Ker\beta_1 &=&\left(1,\frac{1}{2}\right)\oplus \left(\frac{1}{2},1\right),
\\
\Ker\beta_2&=& \left(\frac{3}{2},0\right)\oplus \left(0,\frac{3}{2}\right).
\eea
Indeed, the dimensions of these representations
match the dimensions of the kernel of the operator $\beta$,
\be
\dim\Ker\beta_{0}= 2+2=4,
\quad
\dim\Ker\beta_{1}= 3\cdot 2+2\cdot 3=12,
\quad
\dim\Ker\beta_{2}= 4+4=8.
\ee
Therefore, the dimensions of the representations
of the group $SL(2,\RR)$ are
\bea
\dim V_0 = 8,
\qquad
\dim V_1 = 36,
\qquad
\dim V_2 = 20.
\eea


\section{Geometry }
\setcounter{equation}0
\label{sec:geometry}


Let $M$ be a $n=2\nu$-dimensional oriented compact
Riemannian manifold $M$ without boundary
with a metric $g_{\mu\nu}$. 
Let $\sigma^a$, $a=1,\dots,n$, be an orthonormal frame in $T_x^*M$
and
$\sigma^{ab}=\sigma^a\wedge \sigma^b$, $1\le a<b\le n$,
be the basis in the space of two forms $\Lambda_2$.
We use low case Latin indices from the beginning of the alphabet, $a,b,c,\dots, $
to denote orthonormal frame indices which run over $1,\dots, n$.
We use the Euclidean metric $\delta_{ab}$ to raise these indices.

We define the inner product in the space
of two-forms $\Lambda_2$ by
\be
\la A,B\ra
=\frac{1}{2}A_{ab}B^{ab}.
\ee
We will also use the space $\Lambda_2\otimes \Lambda_2$
of tensors of type $(0,4)$, $A_{abcd}$,
which are anti-symmetric in the first two indices and
anti-symmetric in the last
two indices, with the orthonormal basis
$\sigma^{ab}\otimes\sigma^{cd}$, $1\le a<b\le n$,
$1\le c<d\le n$.
We define the inner product in this
 space 
by
\be
\la A,B\ra
=\frac{1}{4}A_{abcd}B^{abcd}.
\ee

We denote
$g=\det g_{\mu\nu}$ so that
$d\vol= dx g^{1/2}$ is the
Riemannian volume element;
further, we denote
\be
E_{a_1\dots a_n}=g^{1/2}\varepsilon_{a_1\dots a_n},
\ee
where
$\varepsilon_{a_1\dots a_n}$ be the Levi-Civita symbol.
Let $\nabla$ be the Levi-Civita connection
and $R_{abcd}$ be the Riemann curvature tensor;
then Ricci tensor and the scalar curvature are
\be
R_{bd}=R^a{}_{bad },
\qquad
R=R^a{}_a.
\ee
For Einstein manifolds the Ricci tensor
is proprotional to the metric,
\be
R_{\mu\nu}=\lambda g_{\mu\nu}, 
\qquad
R=n\lambda,
\label{150iga}\ee
with some constant $\lambda$.
The Riemann tensor satisfies
the Bianchi identity
\be
\nabla_{[\mu} R^{ab}{}_{\nu\lambda]}=0,
\ee
which gives
\be
\nabla_\alpha R^{\alpha\beta}{}_{\mu\nu}
=2\nabla_{[\mu}R^{\beta}{}_{\nu]},
\qquad
\nabla_\alpha R^\alpha{}_\beta=\frac{1}{2}\nabla_\beta R.
\ee
Here, and every below, the square brackets denote
the anti-symmetrization over all included indices.
This enables one to compute
the Laplacian of the Riemann tensor
\bea
\Delta R^{\alpha\beta}{}_{\mu\nu}
&=&
2\nabla_{[\mu}\nabla^{[\alpha}R^{\beta]}{}_{\nu]}
+2\nabla^{[\alpha}\nabla_{[\mu}R^{\beta]}{}_{\nu]}
-R^{\lambda}{}_{ [\mu} R^{\alpha\beta}{}_{\nu]}{}_{\lambda}
-R^{[\alpha}{}_{\lambda} R^{\beta]\lambda}{}{}_{\mu\nu}
\nn\\
&&
-R^{\alpha\beta}{}_{pq}R^{pq}{}_{\mu\nu}
-4R^{[\alpha}{}_{pq[\mu}R^{\beta] pq}{}_{\nu]}.
\eea


We will study the properties of  various functionals
under the conformal transformation
\be
\hat g_{\mu\nu}=e^{2\omega} g_{\mu\nu}\,,
\ee
where $\omega\in C^\infty(M)$ is a smooth function on $M$.
so that
\be
\hat g^{\mu\nu}=e^{-2\omega}g^{\mu\nu}\,,
\qquad
\hat g=e^{2n\omega}g\,.
\ee
with Christoffel symbols
\be
\hat\Gamma^\alpha{}_{\mu\nu}=\Gamma^\alpha{}_{\mu\nu}
+\delta^\alpha_\mu\omega_\nu
+\delta^\alpha_\nu\omega_\mu
-g_{\mu\nu}\omega^\alpha,
\ee
where $\omega_\mu=\nabla_\mu\omega$.
The curvature tensor transforms by
\bea
\hat R^{\alpha\beta}{}_{\mu\nu}
&=&e^{-2\omega}
\left\{R^{\alpha\beta}{}_{\mu\nu}
-4\delta^{[\alpha}_{[\mu}\omega^{\beta]}{}_{\nu]}
+4\delta^{[\alpha}_{[\mu}\omega^{\beta]}\omega_{\nu]}
-2\delta^{\alpha}_{[\mu} \delta^{\beta}_{\nu]}
\omega_\gamma\omega^\gamma
\right\},
\\
\hat R_{\mu\nu}
&=&R_{\mu\nu}
-g_{\mu\nu}\Delta\omega
-(n-2)\omega_{\mu\nu}
+(n-2)\omega_\mu\omega_\nu
-(n-2)g_{\mu\nu}\omega_\gamma\omega^\gamma,
\\
\hat R
&=&e^{-2\omega} \left\{R
-2(n-1)\Delta\omega
-(n-1)(n-2)\omega_\gamma\omega^\gamma
\right\}
\,.
\eea

For $n>2$ the Weyl tensor is defined by
\be
C^{ab}{}_{cd}=R^{ab}{}_{cd}
-\frac{4}{n-2}\delta^{[a}{}_{[c}R^{b]}{}_{d]}
+\frac{2}{(n-1)(n-2)}
\delta^{[a}{}_{[c}\delta^{b]}{}_{d]}R.
\label{15igax}
\ee
The Weyl tensor transforms
under the conformal transformations
covariantly
\be
\hat C^\alpha{}_{\beta\mu\nu}
= C^\alpha{}_{\beta\mu\nu}.
\ee


By using the canonical decomposition of the Riemann tensor
we get the Bianchi identity for the Weyl tensor
\be
\nabla_{[e}C^{ab}{}_{cd]}
=\frac{4}{n-2} \delta^{[a}{}_{[e}\nabla_{c}R^{b]}{}_{d]}
+\frac{2}{(n-1)(n-2)}
\delta^{[a}{}_{[e}\delta^{b]}{}_{c}\nabla_{d]}R,
\ee
Most of the time we will restrict below to four dimensions, $n=4$;
therefore,
\bea
\Delta C^{\alpha\beta}{}_{\mu\nu}
&=&
2\nabla_{[\mu}\nabla^{[\alpha}R^{\beta]}{}_{\nu]}
+2\nabla^{[\alpha}\nabla_{[\mu}R^{\beta]}{}_{\nu]}
-2\delta^{[\alpha}{}_{[\mu}\Delta R^{\beta]}{}_{\nu]}
+\frac{1}{3}\delta^{[\alpha}{}_{[\mu}\delta^{\beta]}{}_{\nu]}\Delta R
\nn\\
&&
-R^{\lambda}{}_{ [\mu} R^{\alpha\beta}{}_{\nu]}{}_{\lambda}
-R^{[\alpha}{}_{\lambda} R^{\beta]\lambda}{}{}_{\mu\nu}
-R^{\alpha\beta}{}_{pq}R^{pq}{}_{\mu\nu}
-4R^{[\alpha}{}_{pq[\mu}R^{\beta] pq}{}_{\nu]}.
\eea

For Einstein manifolds (\ref{150iga})
this simplifies significantly. The Weyl tensor satisfies
the equations
\be
\nabla_{[e}C^{ab}{}_{cd]}=0,
\ee
and, therefore,
\bea
(-\Delta +2\lambda)C^{\alpha\beta}{}_{\mu\nu}
&=&
C^{\alpha\beta}{}_{pq}C^{pq}{}_{\mu\nu}
+4 C^{[\alpha}{}_{pq[\mu}C^{\beta] pq}{}_{\nu]},
\\
\nabla_{a}C^{a}{}_{bcd}&=&0.
\eea
The solution of this equation significantly depends on the sign of the
Ricci tensor, $\lambda$. For positive $\lambda>0$ there is certain rigidity,
that is, the solutions for small Weyl tensor do not exist. 
For $\lambda<0$ there is much more freedom.
In particular, we have the constraint
\be
||\nabla C||^2
=
\frac{1}{4}\int\limits_M d\vol \left\{
-8\lambda |C|^2
+C^{\mu\nu}{}_{\alpha\beta}C^{\alpha\beta}{}_{pq}C^{pq}{}_{\mu\nu}
+4 C^{\mu\nu}{}_{\alpha\beta}C^{\alpha}{}_{pq\mu}C^{\beta pq}{}_{\nu}
\right\}\ge 0.
\ee


The space of $2$-forms in four dimensions
is canonically decomposed into
the self-dual and anti-self-dual forms
\be
\Lambda_2=\Lambda_2^+\oplus\Lambda^-_2;
\ee
these spaces are three-dimensional, $\dim\Lambda_\pm=3$.

In this case
we can define the
(right)
dual Riemann tensor
\bea
R^*{}_{\mu\nu\alpha\beta}&=&\frac{1}{2}E_{\alpha\beta\rho\sigma}R_{\mu\nu}{}^{\rho\sigma}.
\eea
Notice that the left-dual defined with respect to the first pair of indices
\be
{}^*R_{\mu\nu\alpha\beta}
=\frac{1}{2}E_{\mu\nu\rho\sigma}R^{\rho\sigma}{}_{\alpha\beta},
\ee
is not equal to the right dual.
We also define
the two-sided dual
\bea
{}^*R^*{}^{\mu\nu}{}_{\alpha\beta}&=&\frac{1}{4}E^{\mu\nu\rho\sigma}
E_{\alpha\beta\lambda\kappa}R_{\rho\sigma}{}^{\lambda\kappa}
\nonumber\\
&=&R^{\mu\nu}{}_{\alpha\beta}
-4\delta^{[\mu}{}_{[\alpha}R^{\nu]}{}_{\beta]}
+\delta^{[\mu}{}_{[\alpha}\delta^{\nu]}{}_{\beta]}R.
\eea
Since
$
*^2=1
$
 on $2$-forms, by taking the left
dual of the double dual we obtain
the relation of the right dual and the left dual
\be
R^*_{\mu\nu\alpha\beta}
=
{}^*R_{\mu\nu\alpha\beta}
+2E_{\mu\nu\sigma[\alpha}R^{\sigma}{}_{\beta]}
+\frac{1}{2}E_{\mu\nu\alpha\beta}R.
\ee
From now on we will only use the right dual.


The canonical decomposition of the Riemann tensor
in the irreducible parts in four dimensions, $n=4$, has the form
\bea
R^{ab}{}_{cd}
&=&
C^{ab}{}_{cd}
+2\delta^{[a}{}_{[c}R^{b]}{}_{d]}
-\frac{1}{3}
\delta^{[a}{}_{[c}\delta^{b]}{}_{d]}R.
\label{15iga}
\eea


We also define the dual of the Weyl tensor
\bea
C^*{}_{\mu\nu\alpha\beta}&=&
\frac{1}{2}E_{\rho\sigma\alpha\beta}C_{\mu\nu}{}^{\rho\sigma}.
\eea
It is easy to see that for the Weyl tensor the left and the right dual coincide,
\be
{}^*C=C^*.
\ee
This enables one to define the self-dual and
anti-self-dual parts of the Weyl tensor by
\be
C_{\pm}{}
= \frac{1}{2}\left(C\pm  C^*\right),
\ee
so that
\be
C_{\pm}^* = \pm C_{\pm}.
\ee
Therefore, we have the canonical decomposition of the Weyl tensor into
its self-dual and anti-self-dual parts,
\be
C =C_++C_-,
\qquad
C^* = C_+-C_-.
\ee

The Weyl tensor in four dimensions has ten independent
components.
There are four scalar invariants (invariants under $SO(4)$) of the Weyl tensor,
two quadratic ones,
\be
|C|^2=\frac{1}{4}C_{ab}{}^{cd}C_{cd}{}^{ab}, \qquad
\la C,C^*\ra=\frac{1}{4}C_{ab}{}^{cd}C^{*}_{cd}{}^{ab},
\ee
and two qubic ones
\be
C_{ab}{}^{cd}C_{cd}{}^{ef}C_{ef}{}^{ab},
\qquad
C_{ab}{}^{cd}C_{cd}{}^{ef}C^*_{ef}{}^{ab}.
\ee
The remaining six components are related to the dimension of the
orthogonal group $SO(4)$.

The self-dual and the anti-self-dual parts of the Weyl tensor
have two invariants each,
\be
|C_+|^2=\frac{1}{4}C^+_{ab}{}^{cd}C^+_{cd}{}^{ab}, \qquad C^+_{ab}{}^{cd}C^+_{cd}{}^{ef}C^+_{ef}{}^{ab},
\ee
\be
|C_-|^2=\frac{1}{4}C^-_{ab}{}^{cd}C^{-}_{cd}{}^{ab},\qquad
\qquad
C^-_{ab}{}^{cd}C^-_{cd}{}^{ef}C^-_{ef}{}^{ab}.
\ee
It is easy to see that
\be
\la C_+,C_-\ra=0,
\ee
and, therefore,
\bea
|C^*|^2 &=&|C|^2
 = |C_+|^2+|C_-|^2,
\\
\la C^*,C\ra &=& |C_+|^2-|C_-|^2.
\eea
Then we have for the Riemann tensor
\bea
|\mathrm{Riem}|^2 &=& |C_+|^2+|C_-|^2
+\frac{1}{2}|\Ric|^2
-\frac{1}{12}R^2,
\\
\la {}^*\mathrm{Riem}^*,\mathrm{Riem}\ra
&=& |C_+|^2
+|C_-|^2
-\frac{1}{2}|\Ric|^2
+\frac{1}{6}R^2,
\\
\la\mathrm{Riem}^*,\mathrm{Riem}\ra
&=&
|C_+|^2-|C_-|^2.
\eea


To avoid confusion we use the following
inner product in the space
of {\it symmetric} 2-tensors 
\be
\la A,B\ra=A_{ab}B^{ab}.
\ee


The Weyl tensor acts on the space of two-forms
\be
C: \Lambda_2\to \Lambda_2
\ee
by
\be
(C\varphi)_{\mu\nu}=\frac{1}{2}C_{\mu\nu}{}^{\alpha\beta}\varphi_{\alpha\beta}.
\ee
This operator is symmetric and, therefore, has real eigenvalues;
it can be diagonalized.
This operator commutes with the Hodge star and, therefore,
respects the decomposition of the space of $2$-forms in self-dual
and anti-self-dual forms,
\be
C=C^+\oplus C^-,
\ee
where 
\be
C^\pm: \Lambda_2^\pm\to \Lambda_2^\pm
\ee
 maps self-dual
forms to self-dual forms and anti-self-dual ones to anti-self dual.
Recall that the spaces $\Lambda_2^\pm$ are three-dimensional.
Therefore, the operators $C^\pm$ define two real symmetric
traceless $3\times 3$ matrices. Therefore, they have three real
eigenvalues $W^\pm_i$, $i=1,2,3$, whose sum is equal to zero,
\be
W^\pm_1+W^\pm_{2}+W^\pm_3=0.
\ee
Let $\Sigma^\pm_i$, $i=1,2,3$ be the orthonormal
self-dual and the anti-self-dual
eigen two-forms
\be
C^\pm\Sigma^\pm_i=W^\pm_i \Sigma^\pm_i,
\ee
satisfying
\be
\la\Sigma^+_i,\Sigma^+_j\ra
=\la\Sigma^-_i,\Sigma^-_j\ra=\delta_{ij},
\qquad
\la\Sigma^+_i,\Sigma^-_j\ra=0.
\ee
Then the Weyl tensor has the canonical representation
\be
C_\pm^{abcd}=\sum_{i=1}^3 W^\pm_i \Sigma_i^{ab}\Sigma_i^{cd}
\ee
Therefore, each tensor $C^\pm$
has two real independent invariants, $W^\pm_1$ and $W^\pm_2$, or
\be
I_2=|C_\pm|^2
=\sum_{i=1}^3 (W^\pm_{i})^2
\ee
and
\be
I_3=\sum_{i=1}^3(W^\pm_{i})^3.
\ee


There is a well known relation between the self-dual part of the Weyl
tensor and the scalar curvature in the case of K\"ahler manifolds. 
It turns out that the eigenvalues $W_i^+$ are
\be
W_1^+=\frac{1}{6}R, 
\qquad
W_2^+=-\frac{1}{12}R, 
\qquad
W_3^+=-\frac{1}{12}R,
\ee
so that
\be
C^+_{abcd} = \frac{1}{12}R \Sigma_{abcd},
\ee
where
\be
\Sigma^{abcd}=2\Sigma_1^{ab}\Sigma_1^{cd}
-\Sigma_2^{ab}\Sigma_2^{cd}
-\Sigma_3^{ab}\Sigma_3^{cd},
\ee
and
\be
|C_+|^2=\frac{1}{24}R^2.
\ee
Here $\Sigma_1^+$ is the normalized K\"ahler form and
$\Sigma^+_{2,3}$ are the real and imaginary parts of the
local holomorphic form forming the basis in the orthogonal
complement of K\"ahler form.
In the case of K\"ahler-Einstein metrics with
$R_{ab}=\lambda g_{ab}$ this takes the form
\be
|C_+|^2=\frac{2}{3}\lambda^2.
\ee


Let $B_p(M)$ be the Betti numbers of a $4$-dimensional
manifold $M$.
In the middle dimension, $p=2$, the Betti number $B_2(M)$
can be decomposed into the number of self-dual and anti-self-dual
$2$-forms
\be
B_2(M)=B_2^+(M)+B_2^-(M).
\ee
The two well-known topological invariants of the four-dimensional manifolds
are the Euler characteristic
\be
\chi(M) = 2-2B_1(M)+B_2^+(M)+B_2^-(M)
\ee
and the signature invariant
\be
\tau(M) = B_2^+(M)-B_2^-(M).
\ee
In particular, this means that the invariants
\bea
\chi(M)+\tau(M)
&=&2-2B_1(M)+2B^+_2(M),
\\
\chi(M)-\tau(M)
&=&2-2B_1(M)+2B_2^-(M),
\eea
are even. Therefore, either both $\chi(M)$ and $\tau(M)$
(or $B_2^+(M)$ and $B_2^-(M)$)
are even
or both are odd.

Recall that for simply connected manifolds $B_1(M)=0$. Therefore,
for simply connected manifolds
\be
\chi(M)\ge 2
\ee
and
\bea
B_2^+(M) &=&\frac{1}{2}\left[\chi(M)-2+\tau(M)\right],
\\
B_2^-(M) &=& \frac{1}{2}\left[\chi(M)-2-\tau(M)\right].
\eea
Therefore, for simply connected manifolds
\be
|\tau(M)|\le \chi(M)-2.
\ee



These topological invariants are given by the following integrals
\bea
\chi(M)
&=&
\frac{1}{8\pi^2}\int\limits_M d\vol\;\left\{
|C_+|^2+|C_-|^2
-\frac{1}{2}|\Ric|^2
+\frac{1}{6}R^2\right\},
\\
\tau(M)
&=& \frac{1}{12\pi^2}\int\limits_M d\vol\;\left\{|C_+|^2-|C_-|^2
\right\}.
\eea
Therefore, if $\tau(M)>0$ then $C_+\ne 0$, and
similarly, if $\tau(M)<0$ then $C_-\ne 0$.

Note also the following useful combinations
\bea
2\chi(M)+3\tau(M) &=&
\frac{1}{2\pi^2}\int\limits_M d\vol\;\left\{
|C_+|^2-\frac{1}{4}|\Ric|^2
+\frac{1}{12}R^2\right\},
\\
2\chi(M)-3\tau(M) &=&
\frac{1}{2\pi^2}\int\limits_M d\vol\;\left\{
|C_-|^2-\frac{1}{4}|\Ric|^2
+\frac{1}{12}R^2\right\}.
\eea

Further, for Einstein metrics with the
Ricci tensor $R_{ab}=\lambda g_{ab}$
we have
\bea
\chi(M) &=&
\frac{1}{8\pi^2}\int\limits_M d\vol\;\left\{
|C_+|^2+|C_-|^2\right\}
+\frac{1}{12\pi^2}\lambda^2\vol(M),
\eea
This means, in particular, that a manifold with
negative Euler characteristic, $\chi(M)< 0$,
does not support Einstein metrics.
Furthermore, in this case there hold
\bea
2\chi(M)+3\tau(M) &=&
\frac{1}{2\pi^2}\int\limits_M d\vol\;
|C_+|^2
+\frac{1}{6\pi^2}\lambda^2\vol(M),
\\
2\chi(M)-3\tau(M) &=&
\frac{1}{2\pi^2}\int\limits_M d\vol\;
|C_-|^2
+\frac{1}{6\pi^2}\lambda^2\vol(M),
\eea
which immediately leads to the Hitchin-Thorpe
 inequality for Einstein metrics
\be
\chi(M)\ge \frac{3}{2}|\tau(M)|.
\ee
For the K\"ahler-Einstein metrics this invariant
is determined by the volume
\be
2\chi(M)+3\tau(M) =
\frac{1}{2\pi^2}\lambda^2\vol(M).
\ee


\section{Twisted Exterior Calculus}
\setcounter{equation}0
\label{sec:calculus}


Let $M$ be an even dimensional oriented Riemannian manifold
with a $\Spin^\CC$ structure 
(which always exists in four dimensions).
If $M$ is spin, then there is a globally defined spinor bundle $\cS$,
and we let $L$ be a Hermitian line bundle over $M$.
In general, neither $\cS$ nor $L$ is globally defined, but the 
twisted spinor bundle $\cS\otimes L$ is
(with $L$ playing the role of the square root of the determinant 
line bundle of the $\Spin^\CC$ structure);
all constructions below involve only this bundle.
More generally, we define the bundle
$\cS_p=\cS\otimes\Lambda_p\otimes L$.

Let $A$ be a unitary connection one-form on $L$
and $F=dA$ be its curvature $2$-form.
Let $F_\pm$ be the self-dual and anti-self-dual
parts of the two form $F$; then
\bea
|F|^2 &=& |F_+|^2+|F_-|^2,
\\
\la F^*,F\ra &=&
|F_+|^2-|F_-|^2.
\eea
The Pontryagin invariant is then given by
\be
c_1^2(L)=\frac{1}{4\pi^2}\int\limits_M d\vol
\left(|F_+|^2-|F_-|^2\right).
\label{pontiga}
\ee
We consider the space
of smooth sections $C^\infty(\cS_p)$ and
define the $L^2$ inner product
\be
(\psi,\varphi)=\int\limits_M d\vol \la \psi,\varphi\ra,
\ee
and the corresponding norm $||\psi||^2=(\psi,\psi)$.

Let $\nabla$ be the Levi-Civita connection on the bundle ${\cal S}_p$
of spinor-valued $p$-forms
defined with the connection $A$. The commutator of the
covariant derivatives of spinor valued $p$-forms
is
\bea
[\nabla_\mu,\nabla_\nu]\psi_{\mu_1\dots\mu_p}
=
\left(\frac{1}{4}\cR_{\mu\nu}
+iF_{\mu\nu}
\right)\psi_{\mu_1\dots\mu_p}
-pR_{\mu\nu\alpha[\mu_1}\psi^\alpha{}_{\mu_2\dots\mu_p]},
\eea
where
\be
\cR_{\mu\nu}=R_{ab\mu\nu}\gamma^{ab}.
\ee
In particular, we have  for $p=0$ (pure spinors),
\be
[\nabla_\mu,\nabla_\nu]\psi=\left(\frac{1}{4}\cR_{\mu\nu}
+iF_{\mu\nu}\right)\psi,
\ee
and for $p=1$ (spinor-valued $1$-forms)
\be
[\nabla_\mu,\nabla_\nu]\psi_\alpha=\left(\frac{1}{4}
\cR_{\mu\nu}
+iF_{\mu\nu}\right)\psi_{\alpha}
-R^\beta{}_{\alpha\mu\nu}\psi_\beta.
\ee
Notice that the commutator of covariant derivatives  on the whole 
space $S_*$ introduces a derivative with respect to Grassmann
variables
\bea
[\nabla_\mu,\nabla_\nu]
=\frac{1}{4}\cR_{\mu\nu}
+iF_{\mu\nu}
-R^a{}_{b\mu\nu}\theta^b
\frac{\partial}{\partial \theta^a}.
\eea
This last term acts as a {\it torsion} 
$T^a{}_{\mu\nu}=R^a{}_{b\mu\nu}\theta^b$
on the supermanifold.

The conformal transformation of the orthonormal frame is
\be
\hat e^a_\mu=e^{\omega} e^a_\mu.
\ee
The resulting transformation of the Levi-Civita connection
is
\bea
\hat \nabla_\mu\psi_{\mu_1\dots \mu_p}
&=&\nabla_\mu\psi_{\mu_1\dots\mu_p}
+\left(-\frac{1}{2}\omega^\nu\gamma_{\nu\mu}
-p\omega_\mu\right)
\psi_{\mu_1\mu_2\dots\dots\mu_p}
\nn\\
&&
-p\left(\omega_{[\mu_1}\psi_{|\mu|\mu_2\dots\dots\mu_p]}
-\omega_\nu g_{\mu[\mu_1}
\psi^\nu{}_{\mu_2\dots\dots\mu_p]}\right)
\eea


We introduce the twisted exterior derivative
\be
D_p: {\cal S}_p\to {\cal S}_{p+1}
\ee
(for $0\le p\le n-1$)
acting on spinor-valued forms
by
\be
(D_p\psi)_{\mu_1\dots\mu_{p+1}}
=(p+1)\nabla_{[\mu_1}\psi_{\mu_2\dots\mu_{p+1}]}
\ee
and the twisted coderivative
\be
\bar D_p: {\cal S}_p\to {\cal S}_{p-1}
\ee
(for $1\le p\le n$) by
\be
(\bar D_p\psi)_{\mu_1\dots\mu_{p-1}}=\nabla^\nu\psi_{\nu\mu_1\dots\mu_{p-1}}.
\ee
Formally, we  set
\be
D_n=\bar D_0=0.
\ee
Notice that these operators can be extended to the whole space
$S_*$ by
\be
D=\theta^a\nabla_a,
\qquad
\bar D=\frac{\partial}{\partial \theta^a}\nabla_a.
\ee

The conformal transformation of the twisted exterior derivative is
\bea
(\hat D_p\psi)_{\mu_1\dots\mu_{p+1}}
&=&(D_p\psi)_{\mu_1\dots\mu_{p+1}}
-\frac{1}{2}(p+1)\omega^\nu\gamma_{\nu[\mu_1}
\psi_{\mu_2\dots\dots\mu_{p+1}]},
\\
(\widehat{ {\bar D_p}}\psi)_{\mu_1\dots\mu_{p-1}}
&=& e^{-2\omega}
\left\{({\bar D_p}\psi)_{\mu_1\dots\mu_{p-1}}
+\left[(n-2p)\omega^\mu
-\frac{1}{2}\omega_\nu\gamma^{\nu\mu}\right]
\psi_{\mu\mu_1\dots\dots\mu_{p-1}}
\right\}
\nn\\
\eea


The grading operator $\zeta_p$ anticommutes with the
exterior derivative,
\be
\zeta_p D_p=-D_p\zeta_p,
\qquad
\zeta_p \bar D_p=-\bar D_p\zeta_p,
\ee
The Hodge star intertwines the operators $D$ and $\bar D$,
\be
*_{p-1}\bar D_p =
(-1)^{p+1}
D_{n-p}*_p,
\qquad
*_{p+1}D_p =(-1)^p
\bar D_{n-p}*_p.
\ee
The chirality operator $\Gamma$ commutes with
the exterior derivative
\be
\Gamma D_p=D_p\Gamma,
\qquad
\Gamma \bar D_p=\bar D_p\Gamma.
\ee
Therefore, the operators $D_p$ and $\bar D_p$
act canonically on the subspaces
$\cS_p^+$ and $\cS_p^-$,
and the subspaces $\cS_p^+$ and $\cS_p^-$
are invariant under their action.


An arbitrary natural first order differential operator
\be
\cT_p:\cS_p\to \cS_p
\ee
 has the form
\be
\cT_p
=D_{p-1} B_p+\bar D_{p+1} A_p,
\ee
where 
$B_p: S_p\to S_{p-1}$ and $A_p: S_p\to S_{p+1}$ are
natural endomorphisms; they are linear combinations
of the canonical endomorphisms 
\bea
(B_{p,k}\psi)_{\mu_1\dots\mu_{p-1}}
&=&\gamma_{\nu_1\dots\nu_{k}[\mu_1\dots\mu_{k-1}}
\psi^{\nu_1\dots \nu_{k}}{}_{\mu_{k}\dots \mu_{p-1}]},
\qquad
k=1,\dots, p,
\\
(A_{p,k}\psi)_{\mu_1\dots\mu_{p+1}}
&=&\gamma_{\nu_1\dots\nu_k[\mu_1\dots\mu_{k+1}}
\psi^{\nu_1\dots \nu_k}{}_{\mu_{k+2}\dots\mu_{p+1}]},
\qquad
k=0,1,\dots, p.
\eea
Of particular interest are the
operators $\cT_{p,k}$, $k=0,\dots, p$, defined by
\be
(\cT_{p,k}\psi)_{\mu_1\dots\mu_p}
=
\frac{p!}{k!(p-k)!}
\gamma^{\nu}{}_{\nu_1\dots\nu_k}{}_{[\mu_1\dots\mu_k}
\nabla_{|\nu|}
\psi^{\nu_1\dots \nu_k}{}_{\mu_{k+1}\dots\mu_p]};
\ee
(recall that vertical lines indicate indices excluded from the
antisymmetrization).
Of course these operators can be extended to the whole space $S_*$ by
\be
\cT_k=(-1)^{k(k-1)/2}
\frac{1}{k!}\gamma^{\nu\nu_1\dots\nu_k}{}_{\mu_1\dots\mu_k}
\theta^{\mu_1}\cdots\theta^{\mu_k}
\frac{\partial}{\partial \theta^{\nu_1}}
\cdots\frac{\partial}{\partial \theta^{\nu_k}}
\nabla_{\nu}.
\ee

The  operator $\cD_p=\cT_{p,0}$ is the Dirac operator
\be
(\cD_p\psi)_{\mu_1\dots \mu_p} = \gamma^\nu \nabla_{\nu}\psi_{\mu_1\dots \mu_p},
\ee
and the operator $\cR_p=\cT_{p,1}$ is the Rarita-Schwinger operator
\bea
(\cR_{p}\psi)_{\mu_1\dots\mu_p} &=&
p\gamma^{\nu\lambda}{}_{[\mu_1}\nabla_{|\nu}
\psi_{\lambda|\mu_{2}\dots\mu_p]}.
\eea
By using the identity
\be
\gamma_{\nu\lambda\mu}=
\gamma_\mu\gamma_\nu\gamma_\lambda
-\delta_{\lambda\nu}\gamma_\mu
-\delta_{\mu\nu}\gamma_\lambda
+\delta_{\mu\lambda}\gamma_\nu
\ee
it is easy to show that
the Rarita-Schwinger operator has the form
\bea
\cR_p
&=&
p\cD_p
+
\alpha_{p-1}\cD_{p-1}\beta_p
-D_{p-1}\beta_p
-\alpha_{p-1}\bar D_p.
\eea
Notice that the integer $k$ is also restricted by the dimension,
$
k\le \nu-1.
$
Therefore, in four dimensions, $n=2\nu=4$, the only non-vanishing
operators are the Dirac operator $\cD_p$ and the Rarita-Schwinger
operator $\cR_p$ for $p=1,\dots, 4$.


It is easy to see that
\bea
(\beta_{p+1} D_p\psi)_{\mu_1\dots \mu_p}
&=& (p+1)\gamma^\nu \nabla_{[\nu}\psi_{\mu_1\dots \mu_p]}
\nn\\
&=& \gamma^\nu \nabla_{\nu}\psi_{\mu_1\dots \mu_p}
- p \nabla_{[\mu_1}\gamma^\nu\psi_{|\nu|\mu_2\dots \mu_p]},
\\
(\bar D_{p+1}\alpha_p\psi)_{\mu_1\dots \mu_p}
&=&
(p+1)\nabla^\nu \gamma_{[\nu} \psi_{\mu_1\dots \mu_p]}
\nn\\
&=&
\gamma^\nu\nabla_\nu \psi_{\mu_1\dots \mu_p}
-p \gamma_{[\mu_1} \nabla^\nu\psi_{|\nu|\mu_2\dots \mu_p]},
\eea
Notice that
\bea
(D_{p-1}\beta_p\psi)_{\mu_1\dots \mu_p}
 &=& p\nabla_{[\mu_1}\gamma^\nu\psi_{|\nu|\mu_2\dots \mu_p]},
\\
(\alpha_{p-1}\bar D_p\psi)_{\mu_1\dots \mu_p}
&=& p\gamma_{[\mu_1} \nabla^\nu \psi_{|\nu|\mu_2\dots \mu_p]};
\eea
and, therefore, the Dirac operator can be written in terms of the
twisted exterior derivatives,
\bea
\cD_p &=& \beta_{p+1} D_p+D_{p-1}\beta_p
\nn\\
&=& \alpha_{p-1} \bar D_p+\bar D_{p+1}\alpha_p.
\eea
By using this form and the algebra of the operators $\alpha$ and $\beta$
one can also express the exterior derivatives in terms of the Dirac
operator
\bea
\bar D_p &=& \frac{1}{2}\left(\beta_p \cD_p + \cD_{p-1} \beta_p\right),
\label{fg616}
\\
D_p &=&\frac{1}{2}\left(\alpha_{p} \cD_p + \cD_{p+1} \alpha_p\right).
\label{fg617}
\eea

Under the conformal transformation the Dirac operator transforms
as follows
\bea
(\hat\cD_p\psi)_{\mu_1\dots\mu_p}
&=&e^{-\omega}
\Biggl\{
(\cD_p\psi)_{\mu_1\dots\mu_p}
+\frac{1}{2}(n-2p-1)\omega_\mu \gamma^\mu
\psi_{\mu_1\mu_2\dots\dots\mu_p}
\nn\\
&&
-p\left(\omega_{[\mu_1}
\gamma^\mu\psi_{|\mu|\mu_2\dots\dots\mu_p]}
-\omega_\nu \gamma_{[\mu_1}
\psi^\nu{}_{\mu_2\dots\dots\mu_p]}\right)
\Biggr\}
\eea
This can be written in a more compact form
\bea
(\hat\cD_p e^{-a\omega}\psi)_{\mu_1\dots\mu_p}
&=&e^{-b\omega}
\Biggl\{
(\cD_p\psi)_{\mu_1\dots\mu_p}
-p\left(\omega_{[\mu_1}
\gamma^\mu\psi_{|\mu|\mu_2\dots\dots\mu_p]}
-\omega_\nu \gamma_{[\mu_1}
\psi^\nu{}_{\mu_2\dots\dots\mu_p]}\right)
\Biggr\}
\nn\\
\eea
with 
\be
a=\frac{1}{2}(n-2p-1),
\qquad
b=\frac{1}{2}(n-2p+1).
\ee


The Dirac operator commutes with the Hodge star
\be
\cD_{n-p} *_p = *_{p}\cD_{p}
\ee
and anti-commutes with the chirality operator
\be
\Gamma\cD_p=-\cD_p\Gamma.
\ee
The Dirac operator intertwines the chiral subspaces $\cS_p^+$
and $\cS_p^-$,
\be
\cD^\pm_p: \cS_p^\pm\to \cS_p^\mp
\ee
Therefore, the kernel of the Dirac operator
decomposes accordingly
\be
\Ker\cD_p=\Ker\cD_p^+\oplus \Ker\cD_p^-;
\ee
with a well defined index
\be
\Index\cD^+_p=\dim\Ker\cD_p^+-\dim\Ker\cD_p^-.
\ee

Further, we extend all operators to the total bundle
\be
*, D, \bar D: \cS_*\to \cS_*
\ee
and define the operators
$\cL_{j}: \cS_*\to \cS_*$, $j=1,2$, by
\be
\cL_{1}=D+\bar D,
\qquad
\cL_{2}=i(D-\bar D).
\ee
These operators have the form
\bea
\cL_{1} &=&\left(\frac{\partial}{\partial \theta^a}+\theta^a\right)\nabla_a,
\\
\cL_{2} &=&-i\left(\frac{\partial}{\partial \theta^a}-\theta^a\right)\nabla_a,.
\eea
The  operators $\cL_j$ anti-commute with the
duality involution
\be
\tau_1\cL_1=-\cL_1\tau_1,
\qquad
\tau_2\cL_2=-\cL_2\tau_2,
\ee
where
\be
\tau_1=i^{p(p+1)+\nu}*,
\qquad
\tau_2=i^{p(p-1)+\nu}*,
\ee
 and the grading involution $\zeta$,
\be
\zeta\cL_j=-\cL_j\zeta.
\ee
Notice that the twisted Hodge operators $\cL_j$ commute with the
chirality operator
\be
\Gamma\cL_j=\cL_j\Gamma.
\ee


The restriction of the operators $\cL_j$
to the space of even and odd forms defines the {\it twisted
Hodge operators}
\be
\cL_j^{\rm even}: \cS_{\rm even}\to \cS_{\rm odd},
\qquad
\cL_j^{\rm odd}: \cS_{\rm odd}\to \cS_{\rm even}.
\ee
Similarly,  the restriction of these operators to the
space of self-dual and anti-self-dual forms
defines the {\it twisted signature operator},
\be
\cL_j^{\rm sd}: \cS_{\rm sd}\to \cS_{\rm asd},
\qquad
\cL_j^{\rm asd}: \cS_{\rm asd}\to \cS_{\rm sd}.
\ee
Therefore, these operators intertwine the even and odd forms
and the self-dual and anti-self-dual forms.
The kernels admit two canonical decompositions
\bea
\Ker \cL_j &=& \Ker \cL_j^{\rm even}\oplus \Ker \cL_j^{\rm odd}
\\
&=&\Ker \cL_j^{\rm sd}\oplus \Ker\cL_j^{\rm asd},
\eea
with well defined indices
\bea
\Index \cL_j^{\rm even} &=& \dim\Ker\cL_j^{\rm even}-\dim\Ker\cL_j^{\rm odd},
\\
\Index \cL_j^{\rm sd} &=& \dim\Ker\cL_j^{\rm sd}-\dim\Ker\cL_j^{\rm asd}.
\eea
Notice that there is no corresponding index for the
left and right spinor valued forms.



The square of the operators $D$ and $\bar D$,
more precisely, the products
\bea
Q_p&=& D_{p+1}D_{p}: \cS_p\to \cS_{p+2},
\\
\bar Q_p&=& \bar D_{p-1}\bar D_{p}: \cS_{p}\to \cS_{p-2},
\eea
define the endomorphisms
\bea
(Q_{p}\psi)_{\mu_1\dots\mu_{p+2}} &=&
\frac{1}{2}(p+1)(p+2)
\left(\frac{1}{4}\cR_{[\mu_1\mu_2}
+iF_{[\mu_1\mu_2}\right)
\psi_{\mu_3\dots\mu_{p+2}]},
\nn\\
&&
\hspace{6cm}
0\le p\le n-2,
\label{qpiga}
\\
(\bar Q_{p}\psi)_{\mu_1\dots\mu_{p-2}} &=&
-\frac{1}{2}
\left(
\frac{1}{4}\cR_{\mu\nu}
+iF_{\mu\nu}\right)
\psi^{\mu\nu}{}_{\mu_1\dots\mu_{p-2}},
\qquad 2\le p\le n,
\eea
This defines the operators
\bea
Q 
&=&\frac{1}{2}\left(\frac{1}{4}\cR_{\mu\nu}
+iF_{\mu\nu}
\right)\theta^\mu\theta^\nu,
\\
\bar Q &=&
-\frac{1}{2}
\left(
\frac{1}{4}\cR_{\mu\nu}
+iF_{\mu\nu}\right)\frac{\partial}{\partial \theta^\mu}\frac{\partial}{\partial \theta^\nu}
\eea
Notice that the potentials $Q, \bar Q$
are nilpotent, that is,
\be
\bar Q^j=Q^k=0
\ee
for some integers $j, k$.

The products
\be
D_{p-1}\bar D_{p}: \cS_{p}\to \cS_{p},
\qquad
\bar D_{p+1}D_{p}: \cS_{p}\to \cS_{p},
\ee
are second-order differential operators of the form
\bea
(D_{p-1}\bar D_p\psi)_{\mu_1\dots\mu_{p}} &=&
p\nabla_{[\mu_1}\nabla^\nu\psi_{|\nu|\mu_2\dots\mu_p]},
\\
(\bar D_{p+1}D_p\psi)_{\mu_1\dots\mu_{p}} &=&
\Delta\psi_{\mu_1\dots\mu_{p}}
-p\nabla_{[\mu_1}\nabla^\nu\psi_{|\nu|\mu_2\dots\mu_p]}
\nonumber\\
&&
-p\left(\frac{1}{4}\cR_{\nu[\mu_1}
+i F_{\nu[\mu_1}
+R_{\nu[\mu_1}
\right)\psi^\nu{}_{\mu_2\dots\mu_p]}
\nonumber\\
&&
+\frac{1}{2}p(p-1)R_{\mu\nu[\mu_1\mu_2}
\psi^{\mu\nu}{}_{\mu_3\dots\mu_p]}.
\eea
Obviously, the last two terms are non-zero only for $2\le p\le n$.
The {\it twisted de Rham Laplacian} 
\be
L_p: \cS_{p}\to \cS_{p}
\ee
is
defined by
\be
L_p=-D_{p-1}\bar D_p-\bar D_{p+1}D_p.
\ee
Notice that the Laplacian does not commute with the exterior derivative
and the coderivative.
This operator
 is a Laplace type operator
\be
L_p=-\Delta+W_p,
\label{522iga}
\ee
where $W_p: \cS_p\to \cS_p$ is the Weitzenb\"ock endomorphism,
defined by
\bea
(W_p\psi)_{\mu_1\dots\mu_{p}} &=&
p\left(\frac{1}{4}\cR_{\nu[\mu_1}
+R_{\nu[\mu_1}
+i F_{\nu[\mu_1}\right)\psi^\nu{}_{\mu_2\dots\mu_p]}
\nonumber\\
&&
-\frac{1}{2}p(p-1)R_{\mu\nu[\mu_1\mu_2}\psi^{\mu\nu}{}_{\mu_3\dots\mu_p]},
\qquad
2\le p\le n.
\eea
For $p=0$ the endomorphism $W_0$ vanishes, $W_0=0$, and for
$p=1$ it has the form
\bea
(W_1\psi)_{\mu} &=&
\left(\frac{1}{4}\cR_{\nu\mu}+R_{\nu\mu}+i F_{\nu\mu}\right)\psi^\nu.
\eea
It is easy to see that the Weitzenb\"ock term on the whole space $S_*$
has the form
\bea
W &=& \left(\frac{1}{4}\cR_{\mu a}
+R_{a\mu}
+i F_{\mu a}\right)
\theta^a\frac{\partial}{\partial \theta^\mu}
-\frac{1}{2}R_{\mu\nu ab}
\theta^a\theta^b\frac{\partial}{\partial \theta^\nu}
\frac{\partial}{\partial \theta^\mu},
\eea


It is easy to see that the Laplacian commutes
with the chirality operator,
with the Hodge star
and the grading operator,
\be
\Gamma L_p=L_p\Gamma,
\qquad
L_{n-p}*_p=*_pL_p,
\qquad
\zeta_p L_p=L_p\zeta_p.
\ee


The  square of the Dirac operator is
a Laplace type operator,
\be
\cD^2_p =\Delta-V_p,
\label{620iga}
\ee
where $V=-\gamma^{\mu\nu}\nabla_\mu\nabla_\nu$ is an endomorphism defined by
\bea
(V_p\psi)_{\mu_1\dots\mu_p}
&=&\left(\frac{1}{4}R
-\frac{i}{2}\cF\right)\psi_{\mu_1\dots\mu_p}
+\frac{p}{2}\cR_{\alpha[\mu_1}\psi^\alpha{}_{\mu_2\dots\mu_p]};
\eea
here
\be
\cF=F_{\mu\nu}\gamma^{\mu\nu},
\ee
In particular,  for $p=0$ the last term vanishes, so,
$V_0$ reads
\bea
V_0\psi
&=&
\left(\frac{1}{4}R-\frac{1}{2}i\cF\right)
\psi
\eea
and for $p=1$, $V_1$ it has the form
\bea
(V_1\psi)_{\mu}
&=&
\left(\frac{1}{4}R-\frac{1}{2}i\cF\right)
\psi_{\mu}
+\frac{1}{2}\cR_{\nu\mu}\psi^\nu.
\eea
This potential term can be written in the form
\be
V=\frac{1}{4}R
-\frac{i}{2}\cF
+\frac{1}{2}\cR_{a b}
\theta^b\frac{\partial}{\partial \theta^a}.
\ee

The squares of the operators $\cL_j$ define the {\it twisted
Hodge Laplacians}
\be
\cL^2_j =
\Delta-\cW_j,
\ee
where
\bea
\cW_{1,2} &=& W\mp(Q+\bar Q).
\eea


\section{Global Analysis}
\setcounter{equation}0
\label{sec:analysis}

It is easy to see that the coderivative $\bar D$ is the negative formal adjoint of $D$, 
\be
\bar D_p=-D^*_{p-1},
\ee
and the Dirac operator is anti-self-adjoint
\be
\cD_p^*=-\cD_p.
\ee
The Laplacian, $L_p$, and the square of the 
Dirac operator, $-\cD_p^2$, are elliptic self-adjoint
second-order
partial differential operators of Laplace type acting on
smooth sections of vector bundles over a
compact manifold.

The twisted de Rham Laplacians are non-negative operators
since
\bea
(\psi, L_p\psi)&=&
||D_p\psi||^2+||\bar D_p\psi||^2.
\label{51igax}
\eea
Therefore, a spinor valued $p$-form $\psi$ is harmonic,
\be
L_p\psi=0
\ee
if and only if
\be
D_p\psi=\bar D_p\psi=0.
\ee
By using (\ref{522iga}) we have
\bea
(\psi, L_p\psi)
&=&\int\limits_M d\vol
\left\{|\nabla\psi|^2
+\la \psi, W_p\psi\ra
\right\},
\label{51iga}
\eea
where for
$2\le p\le n$,
\bea
\la \psi, W_p\psi\ra&=&
R_{\alpha\beta}G^{\alpha\beta}
-iF_{\alpha\beta}\Omega^{\alpha\beta}
-\frac{1}{4}R_{\alpha\beta\mu\nu}
\left(2J^{\alpha\beta\mu\nu}
+A^{\alpha\beta\mu\nu}
\right),
\eea
where
\bea
G^{\alpha\beta}&=&
\frac{1}{(p-1)!}
\bar\psi^{(\alpha|\mu_1\dots\mu_{p-1}|}
\psi^{\beta)}{}_{\mu_1\dots\mu_{p-1}},
\label{98igax}
\\
\Omega^{\alpha\beta}&=&
\frac{1}{(p-1)!}
\bar\psi^{[\alpha|\mu_1\dots\mu_{p-1}|}
\psi^{\beta]}{}_{\mu_1\dots\mu_{p-1}},
\label{98igaxz}\\
J^{\alpha\beta\mu\nu}&=&\frac{1}{(p-2)!}
\bar\psi^{\alpha\beta\mu_1\dots\mu_{p-2}}
\psi^{\mu\nu}{}_{\mu_1\dots\mu_{p-2}},
\label{99igax}\\
A^{\alpha\beta\mu\nu}&=&\frac{1}{(p-1)!}
\bar\psi^{\alpha\mu_1\dots\mu_{p-1}}\gamma^{\mu\nu}
\psi^\beta{}_{\mu_1\dots\mu_{p-1}}.
\label{910igax}
\eea
Of course, the tensors $G_p$, $\Omega_p$ and $A_p$ are only defined for $p\ge 1$
and the tensor $J_p$ is only defined for $p\ge 2$.
So, we define the tensors $G_p$ and $\Omega_p$ to be zero for $p=0$
and the tensor $J_p$ to be zero for  $p=0,1$.
These tensors have the following symmetries
\bea
G_{\alpha\beta} &=& G_{\beta\alpha}= \overline{G_{\alpha\beta}},
\\
\Omega_{\alpha\beta} &=& -\Omega_{\beta\alpha}
=-\overline{\Omega_{\alpha\beta}},
\\
J_{\alpha\beta\mu\nu}&=&-J_{\beta\alpha\mu\nu}
=-J_{\alpha\beta\nu\mu}
=\overline{J_{\mu\nu\alpha\beta}},
\\
A_{\alpha\beta\mu\nu}&=&
-A_{\alpha\beta\nu\mu}
=-\overline{A_{\beta\alpha\mu\nu}}.
\eea
Therefore, the tensors $G_{\mu\nu}$ and
$A_{[\alpha\beta]\mu\nu} $ are real and the tensors 
$\Omega_{\alpha\beta}$ 
and $A_{(\alpha\beta)\mu\nu} $ are imaginary.
Notice also that
\bea
J^\mu{}_{\beta\mu\nu} &=& (p-1)\left(G_{\beta\nu}+\Omega_{\beta\nu}\right),
\\
J^{\mu\nu}{}_{\mu\nu} &=& (p-1)G^\mu{}_{\mu}
=p(p-1)|\psi|^2.
\eea

We will find it useful to rewrite it in terms
of the Weyl tensor
\bea
&&
\la \psi, W_p\psi\ra
=
\frac{1}{12}R
\left(2p(p-1)|\psi|^2
+A^{\mu\nu}{}_{\mu\nu}
\right)
-iF_{\alpha\beta}\Omega^{\alpha\beta}
\nn\\
&&\qquad
+R_{\alpha\beta}\left((2-p)G^{\alpha\beta}
-\frac{1}{2}H^{\alpha\beta}
\right)
-\frac{1}{4}C_{\alpha\beta\mu\nu}
\left(2J^{\alpha\beta\mu\nu}
+A^{\alpha\beta\mu\nu}
\right),
\label{818igax}
\eea
where
\be
H_{\alpha\beta}=\frac{1}{2}
\left(A_{\alpha}{}_\mu{}_{\beta}{}^\mu
-A_\mu{}_{\alpha\beta}{}^\mu
\right).
\label{917igax}
\ee
So, for $p=0$ this endomorphism vanishes,
$W_0=0$,
and for $p=1$ it has the form
\bea
&&\la \psi, W_1\psi\ra
=
\frac{1}{12}R
\bar\psi^{\mu}\gamma_{\mu\nu}\psi^\nu{}
-iF_{\alpha\beta}\bar\psi^{\alpha}\psi^\beta{}
\\
&&\qquad
+R_{\alpha\beta}\left(
\bar\psi^{\alpha}\psi^\beta{}
-\frac{1}{4}
\bar\psi_{\mu}\gamma^{\mu\alpha}\psi^\beta{}
-\frac{1}{4}
\bar\psi^{\beta}\gamma^{\alpha\mu}\psi_\mu{}
\right)
-\frac{1}{4}C_{\alpha\beta\mu\nu}
\bar\psi^{\alpha}\gamma^{\mu\nu}
\psi^\beta{}.
\nn\eea

Therefore, if
\be
\int\limits_M d\vol \la \psi, W_p\psi\ra >0,
\ee
in particular, if
 $\la \psi, W_p\psi\ra>0$
 everywhere, then there are no harmonic forms,
$\Ker L_p=\{0\}$, and the operator $L_p$ is strictly
positive, $L_p>0$.

Since both Laplacians commute with the Hodge star the Weitzenb\"ock
endomorphism also commutes with it,
\be
W_{n-p}*_p=*_pW_p,
\ee
that is,
\be
\la *\psi, W *\psi\ra=\la\psi,W\psi\ra.
\ee


Similarly, the twisted Hodge operators
$\cL_{j}$
are anti-self-adjoint,
\be
\cL_j^*=-\cL_j.
\ee
and, therefore,
for the twisted Hodge Laplacians
we have
\be
||\cL_{1,2}\psi||^2 
=\sum_{p=0}^n||\nabla\psi_p||^2
+(\psi, \cW_j\psi),
\ee
where
\bea
(\psi, \cW_{1,2}\psi)
&=&
\sum_{p=1}^n\left( \psi_p, W_p\psi_p\right)
\nn\\
&&
\mp\left(\sum_{p=0}^{n-2}
\left(\psi_{p+2}, Q_p\psi_p\right)
+\sum_{p=2}^n
\left(\psi_{p-2}, \bar Q_p\psi_p\right)
\right).
\label{57igax}
\eea
In particular, this means for any $\psi\in S_*$,
\be
\sum_{p=0}^n||\nabla\psi_p||^2
+\sum_{p=1}^n(\psi_p, W_p\psi_p)
\ge
\left|\sum_{p=0}^{n-2}(\psi_{p+2}, Q_p\psi_p)
+\sum_{p=2}^n(\psi_{p-2}, \bar Q_p\psi_p)
\Biggr\}
\right|.
\ee


The negative square of the Dirac operator is
a non-negative operator
since
\bea
-(\psi,\cD^2_p\psi) &=& ||\cD_p\psi||^2.
\label{532igax}
\eea
By using (\ref{620iga}) we have
\bea
||\cD_p\psi||^2
&=&
\int\limits_M d\vol\left\{
|\nabla\psi|^2
+\la \psi, V_p\psi\ra
\right\},
\label{532iga}
\eea
where for $1\le p\le n$,
\bea
\la \psi, V_p\psi\ra &=&\frac{1}{4}R|\psi|^2
-\frac{i}{2}F_{\mu\nu}B^{\mu\nu}
-\frac{1}{2}R_{\alpha\beta\mu\nu}
A^{\alpha\beta\mu\nu},
\eea
where
\be
B^{\mu\nu}=\frac{1}{p!}
\bar\psi^{\mu_1\dots\mu_{p}}\gamma^{\mu\nu}
\psi_{\mu_1\dots\mu_{p}}.
\label{929igax}
\ee
Notice that 
\be
A^\alpha{}_{\alpha\mu\nu}= pB_{\mu\nu}.
\ee

We rewrite this in terms of the Weyl tensor
\be
\la \psi, V_p\psi\ra 
=R\left(\frac{1}{4}|\psi|^2
+\frac{1}{6}A^{\mu\nu}{}_{\mu\nu}\right)
-\frac{i}{2}F_{\mu\nu}B^{\mu\nu}
-R_{\alpha\beta}H^{\alpha\beta}
-\frac{1}{2}C_{\alpha\beta\mu\nu}
A^{\alpha\beta\mu\nu}.
\label{925iga}
\ee

For $p=0$ this takes the form
\bea
\la \psi, V_0\psi\ra &=&
\frac{1}{4}R|\psi|^2
-\frac{1}{2}iF_{\mu\nu}\bar\psi^{}
\gamma^{\mu\nu}\psi_{},
\label{59iga}
\eea
and for $p=1$ this reads
\bea
\la \psi, V_1\psi\ra &=&R\left(\frac{1}{4}|\psi|^2
+\frac{1}{6}
\bar\psi^{\mu}\gamma_{\mu\nu}\psi^\nu{}
\right)
-\frac{1}{2}iF_{\mu\nu}\bar\psi^{\alpha}
\gamma^{\mu\nu}\psi_{\alpha}
\nn\\
&&
-\frac{1}{2}R_{\alpha\beta}\left(
\bar\psi_{\mu}\gamma^{\mu\alpha}\psi^\beta{}
+\bar\psi^{\beta}\gamma^{\alpha\mu}\psi_\mu{}
\right)
-\frac{1}{2}C_{\alpha\beta\mu\nu}
A^{\alpha\beta\mu\nu}.
\eea

Therefore,
if
\be
\int\limits_M d\vol \la \psi, V_p\psi\ra > 0,
\ee
in particular, if $\la \psi, V_p\psi\ra>0$
everywhere, then
the Dirac operator does not have zero modes,
$\Ker \cD_p=\{0\}$, and the operator $-\cD^2_p$ is strictly
positive,
$
-\cD^2_p>0.
$

The Dirac operator also
 commutes with the Hodge star, therefore, so does
the Weitzenb\"ock endomrohism $V$,
\be
V_{n-p}*_p=*_pV_p,
\ee
that is,
\be
\la *\psi, V *\psi\ra=\la\psi, V\psi\ra.
\ee


The spectra of Laplace type operators on compact manifolds
are
real discrete non-negative with finite multiplicities
with well defined heat kernels. We will be interested in the following heat traces
\bea
X_p(t) &=& \Tr\exp(-tL_p),
\\
I_p(t)&=& \Tr\exp(t\cD_p^2),
\\
J_p(t)&=& \Tr\Gamma\exp(t\cD_p^2)
=\Tr\exp(t\cD^2_{p,+})
-\Tr\exp(t\cD^2_{p,-}),
\eea
We can also define
the {\it total heat traces} of the operators $\cL_1^2$ and $\cL_2^2$,
\bea
\Theta_j(t)&=& \Tr\exp(t\cL_j^2),
\\
\Psi_j(t) &=&\Tr\zeta\exp(t\cL_j^2)
=\Tr\exp(t\cL_{j,\rm even}^2)-\Tr\exp(t\cL_{j,\rm odd}^2),
\\
\Phi_j(t)&=&\Tr\tau_j\exp(t\cL_j^2)
=\Tr\exp(t\cL_{j,\rm sd }^2)
-\Tr\exp(t\cL_{j,\rm asd}^2).
\eea


Let $A$ be an anti-self-adjoint operator with a nonpositive square
$A^2\le 0$ and
$C$ be a self-adoint involution anticommuting with the operator $A$,
that is,
\be
AC=-CA.
\ee
Then it is easy to see that
\be
C A^2\exp(tA^2)
= - ACA\exp(tA^2)
= - A^2C\exp(tA^2)
\ee
and, therefore,
\be
\frac{d}{dt}\Tr C\exp(tA^2)=-\Tr CA^2\exp(tA^2)=0,
\ee
that is, the invariant
\be
\Index A=\Tr C\exp(tA^2)
\ee
is constant.

Since the Dirac operator anticommutes with the chirality operator
$\Gamma$, this means that
the invariant $J_p$ is equal to the
index of the Dirac operator, $\cD^+_p$,
\be
J_p=\Index \cD^+_p.
\ee
Similarly, since the operators $\cL_j$ anticommute with the
respective duality operator $\tau_j$,
the invariants
$\Phi_j$, $j=1,2,$ are the indices of the {\it twisted signature operators},
\be
\Phi_j=\Index \cL^{\rm sd}_{j},
\ee
and since the operators $\cL_j$ anti-commute with the
grading
involution $\zeta$,
the invariants $\Psi_j$, $j=1,2$, are the indices of the {\it twisted
Hodge operators},
\be
\Psi_j=\Index \cL^{\rm even}_j.
\ee
Note that, since the index of the Euler characteristic operator 
twisted by a vector bundle $E$ equals $\operatorname{rank}(E)\,\chi(M)$,
these indices are determined by the Euler characteristic,
$
\Psi_1=\Psi_2
=2^\nu\chi(M).
$

The asymptotic expansion of the heat kernel diagonal of a Laplace type operator
$L$ on a closed manifold as $t\to 0$ is
\be
U^{\rm diag}_L(t)\sim (4\pi t)^{-n/2}
\sum_{k=0}^\infty \frac{(-1)^k}{k!}t^k a_k(L),
\ee
where $a_k$ are the local heat kernel coefficients.
The asymptotic expansion of the heat trace with an endomorphism
$E$ is
\be
\Tr E\exp(-tL)\sim (4\pi t)^{-n/2}
\sum_{k=0}^\infty \frac{(-1)^k}{k!}t^k A_k(L,E),
\ee
where
\be
A_k(L,E)=\int_M d\vol \tr \left[E a_k(L)\right].
\ee


Therefore, the indices of the above operators are
non-zero only in even dimension $n=2\nu$,
\bea
\Index \cD^+_p &=& (-1)^{\nu}\frac{1}{(4\pi)^{\nu}\nu!} A_{\nu}(-\cD^2_p,\Gamma),
\\
\Index \cL^{\rm even}_j 
&=& (-1)^{\nu}\frac{1}{(4\pi)^{\nu}\nu!} A_{\nu}(-\cL^2_j,\zeta),
\\
\Index \cL^{\rm sd}_j 
&=& (-1)^{\nu}\frac{1}{(4\pi)^{\nu}\nu!} A_{\nu}(-\cL^2_j,\tau_j).
\eea
In four dimensions these are determined by the coefficient $A_2$.


\section{Monopoles}
\setcounter{equation}0
\label{sec:monopoles}


We restrict ourselves in this section to four dimensions, $n=4$,
and to pure right spinors, $\psi_+\in S_0^+$, 
that is, we consider the case
$p=0$.
Let $F_+$ be the self-dual part of the
curvature of the $U(1)$ connection.
We define, following \cite{witten1994},
the Seiberg-Witten functional
\bea
S_{SW}(g,\psi_+,A) &=& ||\cD\psi_+||^2
+\frac{1}{2}||F_+ + iB_+||^2.
\eea
We notice the relation (\ref{450uni}) (valid in $n=4$)
\be
|iB_+|^2
=2|\psi_+|^4,
\ee
By using eq. (\ref{532iga}) and (\ref{59iga})
we compute
\bea
-\la \psi_+,\cD^2\psi_+\ra 
&=&
-\la\psi_+,\Delta \psi_+\ra
+\frac{1}{4}R|\psi_+|^2
-i\la F_+,B_+\ra
\nn\\
&=&
-\frac{1}{2}\Delta|\psi_+|^2
+|\nabla\psi_+|^2
+\frac{1}{4}R|\psi_+|^2
-i\la F_+,B_+\ra
\eea
Therefore, we have
\bea
||\cD\psi_+||^2
&=& \int\limits_M d\vol
\left\{|\nabla\psi_+|^2
+\frac{1}{4}R|\psi_+|^2
-i\la F_+,B_+\ra
\right\}
\eea
and
\bea
S_{SW}(g,\psi_+,A) &=&
\int\limits_M d\vol
\left\{|\nabla\psi_+|^2
+\frac{1}{4}R|\psi_+|^2
+|\psi_+|^4
+\frac{1}{2}|F_+|^2
\right\}.
\eea
The  functional $S(g,\psi_+,A)$
is nonnegative and has absolute minima
equal to zero
when the fields satisfy the {\it monopole equations}
\bea
\cD\psi_+ &=& 0,
\\
F_+ &=& -iB_+.
\eea
Notice that if the scalar curvature is positive,
$R>0$,
then the solution is
trivial
\be
\psi_+=F_+=0,
\ee
which defines anti-self-dual $U(1)$ connections.

It is easy to see that the solution  $\psi_+$ 
of the monopole equations
is bounded
by the scalar curvature. Indeed, let 
$(\psi_+, A)$ be a solution. Then, obviously, 
$S(g,\psi, A)=\cD\psi_+=0$ and, therefore,
\bea
|\nabla\psi_+|^2
+\frac{1}{4}R|\psi_+|^2
+2|\psi_+|^4=\frac{1}{2}\Delta|\psi_+|^2.
\eea
Suppose that $|\psi_+(x)|^2$ attains a maximum at a point
$x_0\in M$. Then at this point
$\Delta|\psi|^2\le 0$; therefore,
\be
|\psi_+(x_0)|^2\le -\frac{1}{8}R(x_0).
\ee
Therefore, for any $x\in M$
\be
|\psi_+(x)|^2\le -\frac{1}{8}R_{\rm min},
\ee
where
\be
R_{\rm min}=\min_{x\in M} R(x);
\ee
(recall that $R_{\rm min}<0$.)

This functional is invariant under the diffeomorphisms
$\varphi: M\to M$,
\be
g\mapsto \varphi^* g,
\qquad
\psi_+\mapsto \varphi^*\psi_+,
\qquad
A\mapsto \varphi^* A,
\ee
where $\varphi^*$ is the pullback of the diffeomorphism
$\varphi$,
and the $U(1)$ gauge transformations
\be
\psi_+\mapsto \exp(if)\psi_+,
\qquad
A\mapsto A+df,
\ee
where $f\in C^\infty(M)$ is a smooth function.
These invariance transformations form a group $\cG$
called the {\it gauge group}.
The Seiberg-Witten {\it moduli space} $\cM$ is defined as the space of solutions
of these equations modulo the gauge group $\cG$.

By using the
decomposition
\be
|F|^2=|F_+|^2+|F_-|^2,
\ee
this functional
can be rewritten in the form
\be
S_{SW}(g,\psi_+,A)
=
\int\limits_M d\vol
\left\{|\nabla\psi_+|^2
+\frac{1}{4}R|\psi_+|^2
+|\psi_+|^4
+\frac{1}{4}|F|^2\right\}
+\pi^2 c_1^2(L),
\ee
where
$
c_1^2(L)=p_1(L)
$
is the Pontryagin invariant (\ref{pontiga}).
Therefore, the solutions of the monopole equations
are related to the Pontryagin invariant by
\be
c_1^2(L)
=-\frac{1}{\pi^2}
\int\limits_M d\vol
\left\{|\nabla\psi_+|^2
+\frac{1}{4}R|\psi_+|^2
+|\psi_+|^4
+\frac{1}{4}|F|^2\right\}.
\ee

It has been shown that the moduli space for the Seiberg–Witten equations
generically has the dimension
\cite{witten1994}
\bea
\dim {\cal M} &=&c_1^2(L)-\frac{1}{4}\left(2\chi(M)+3\tau(M)\right)
\label{813iga}
\\
&=&
\frac{1}{8\pi^2}\int\limits_M d\vol\;
\left\{2\left(|F_+|^2-|F_-|^2\right)
-|C_+|^2+\frac{1}{4}|\Ric|^2
-\frac{1}{12}R^2
\right\}.
\nn
\eea

\bea
2\chi(M)+3\tau(M) &=&
\frac{1}{2\pi^2}\int\limits_M d\vol\;\left\{
|C_+|^2-\frac{1}{4}|\Ric|^2
+\frac{1}{12}R^2\right\},
\eea
It vanishes when
\be
c_1^2(L)=\frac{1}{4}\left(2\chi(M)+3\tau(M)\right),
\ee
that is, when
\be
\int\limits_M d\vol\;
\left(
|F_+|^2-|F_-|^2
\right)
=\frac{1}{2}\int\limits_M d\vol\;\left\{
|C_+|^2
-\frac{1}{4}|\Ric|^2
+\frac{1}{12}R^2\right\}.
\ee
In this case the moduli space is just a collection of finite
number of points, $P_i$, $i=1,\dots, N$.
Each point  has a sign $\varepsilon(P_i)$ associated with it, coming from the sign of
the determinant of the elliptic operator whose index gave the
dimension of the moduli space. The sum of these signs
\be
W=\sum_{i=1}^N\varepsilon(P_i).
\ee
is a
topological invariant under some conditions (like $B_{2,+}(M)\ge 2$)
\cite{witten1994}.

%
%
%


We propose to generalize the Seiberg-Witten monopoles in
the following way.
We consider the Dirac operator acting on the
right spinor valued $p$-forms in four dimensions, $n=4$,
and $\nu=2$. 
We use 
(\ref{532iga}) 
to define a functional
\bea
&&
S_D(\psi,g,A)
=
||\cD\psi_+||^2+\frac{1}{2}||F_++iB_+||^2
\\
&&\qquad
=
\int\limits_M d\vol\Biggl\{
|\nabla\psi_+|^2
+\frac{1}{2}|B_+|^2
+R\left(\frac{1}{4}|\psi_+|^2
+\frac{1}{6}A_+^{\mu\nu}{}_{\mu\nu}\right)
\nn\\
&&\qquad
-R_{\alpha\beta}H_+^{\alpha\beta}
-\frac{1}{2}C_{\alpha\beta\mu\nu}
A_+^{\alpha\beta\mu\nu}
+\frac{1}{2}|F_+|^2
\Biggr\};
\nn\eea
here the tensors $A_{\mu\nu\alpha\beta}$,
$B_{\mu\nu}$, and $H_{\mu\nu}$
are defined by
(\ref{910igax}), (\ref{929igax}), 
and (\ref{917igax}).
Notice that,
by using (\ref{2103iga}), or, more specifically, (\ref{2122iga}), we
 compute
\bea
|B_+|^2 &=&
\frac{4}{(p!)^2}
(\bar\psi_+^{\mu_1\dots\mu_{p}}{}
\psi^+_{\nu_1\dots\nu_{p}})
(\bar\psi_+^{\nu_1\dots\nu_{p}}
\psi^+_{\mu_1\dots\mu_{p}})
-2|\psi_+|^4
\nn\\
&=&
4\tr(\Phi^2)-2(\tr \Phi)^2,
\eea
where
\be
\Phi=\frac{1}{p!}\psi^+_{\nu_1\dots\nu_{p}}
\otimes
\bar\psi_+^{\nu_1\dots\nu_{p}}.
\ee
This functional is non-negative; it vanishes when
the fields satisfy the monopole equations
\bea
\cD\psi_+&=&0,
\\
F_+ &=& -iB_+.
\eea
It is important to stress that the quantity $|B_+|^2$
vanishes whenever $\Phi$ is proportional to the unit matrix.
Therefore, unlike in the Seiberg-Witten case there is no
coercive term $|\psi|^4$ that control compactness of the moduli 
space.


Next, we consider the twisted de Rham Laplacian acting on the
right spinor valued forms.
We use the equations (\ref{51iga}) and (\ref{818igax})
to define the functional
\bea
&&
S_{dR}(g,\psi,A) =
||D_p\psi_+||^2+||\bar D_p\psi_+||^2
+||F+i\Omega||^2
\\
&&\qquad
=
\int\limits_M d\vol
\Biggl\{
|\nabla\psi_+|^2
+|\Omega|^2
+\frac{1}{12}R
\left(2p(p-1)|\psi_+|^2
+A^{\mu\nu}{}_{\mu\nu}
\right)
\nn\\
&&\qquad
+R_{\alpha\beta}\left((2-p)G^{\alpha\beta}
-\frac{1}{2}H^{\alpha\beta}
\right)
-\frac{1}{4}C_{\alpha\beta\mu\nu}
\left(2J^{\alpha\beta\mu\nu}
+A^{\alpha\beta\mu\nu}
\right)
+|F|^2
\Biggr\};
\nn
\eea
where the tensors $G_{\mu\nu}$, $\Omega_{\mu\nu}$, 
$J_{\mu\nu\alpha\beta}$, $A_{\mu\nu\alpha\beta}$ and $H_{\mu\nu}$ 
are defined by (\ref{98igax})--(\ref{910igax})
and 
(\ref{917igax}). 
The zeros of this functional are realized by monopole
equations
\bea
D_p\psi_+ &=&0,
\\
\bar D_p\psi_+ &=& 0,
\\
F &=& -i\Omega.
\eea
It is easy to see that these equations are overdetermined;
they lead to the integrability constraints
\bea
Q_{p}\psi_+&=& 0,
\\
\bar Q_{p}\psi_+&=&0,
\\
d\Omega &=& 0.
\eea
where $Q_p$ and $\bar Q_p$ are endomorphisms
defined by (\ref{qpiga}).

We believe it might be worthwhile investigating these
generalized monopole equations and their consequences.
In particular, on special manifolds when there is some
control on the Weitzenb\"ock endomorphisms in these
formulas like Einstein metric or K\"ahler manifolds.


\end{document}